\documentclass[11pt]{amsart}

\usepackage[english]{babel}

\usepackage{amsmath, amsthm, amssymb, amsfonts, amsopn}
\usepackage[foot]{amsaddr}

\usepackage{mathrsfs, mathtools, dsfont, esint, bm}

\usepackage{array, multirow, makecell, booktabs}
\setcellgapes{3pt} % vertical spacing in tables (optional)
\renewcommand\arraystretch{1.3}

\usepackage{algorithm2e}

\usepackage{graphicx, import}
\usepackage{pgf, tikz}
\usetikzlibrary{arrows}

\usepackage[skip=2pt,font=scriptsize]{caption}
\usepackage{subcaption}
\usepackage{floatrow}

\usepackage{enumerate}
\usepackage{enumitem}

\usepackage{xcolor}

\usepackage{hyperref}
\usepackage{cleveref}
\hypersetup{
    colorlinks,
    linkcolor={red!80!black},
    citecolor={green!50!black}
}

\usepackage[margin=1in]{geometry}

\theoremstyle{plain}
\newtheorem{theorem}{Theorem}[section]

\newtheorem{proposition}[theorem]{Proposition}

\theoremstyle{definition}
\newtheorem{definition}[theorem]{Definition}

\numberwithin{equation}{section}
\numberwithin{figure}{section}
\numberwithin{table}{section}

\newcommand{\R}{\mathbb{R}}

\newcommand{\D}{\mathcal{D}}

\newcommand{\LL}{\mathscr{L}}

\renewcommand{\Sigma}{\boldsymbol{\sigma}}
\renewcommand{\L}{\mathsf{L}}

\renewcommand{\u}{\boldsymbol{u}}

\DeclareMathOperator{\argmin}{argmin}

\newcommand{\x}{\boldsymbol{x}}

\title[HYCO: Hybrid-Cooperative Learning]{HYCO: Hybrid-Cooperative Learning \\ for Data-Driven PDE Modeling}
\author[L. Liverani, E. Zuazua]{Lorenzo Liverani$^\dagger$}
\address{$^\dagger$Chair for Dynamics, Control, Machine Learning, and Numerics (Alexander von Humboldt Professorship), Department of Mathematics, Friedrich-Alexander-Universit\"at Erlangen-N\"urnberg, 91058 Erlangen, Germany.}

\author{Enrique Zuazua$^\dagger$$^*$$^\S$}
\address{$^*$Departamento de Matem\'aticas, Universidad Aut\'onoma de Madrid, 28049 Madrid, Spain.}
\address{$^\S$Chair of Computational Mathematics, Universidad de Deusto. Av. de las Universidades, 24, 48007 Bilbao, Basque Country, Spain.}
\email{lorenzo.liverani@fau.de}
\email{enrique.zuazua@fau.de}

\makeatletter
\@ifundefined{subjclassname@2020}{%
  \expandafter\gdef\csname subjclassname@2020\endcsname{%
    \textup{2020} Mathematics Subject Classification}%
}{}
\makeatother

\begin{document}

% REQUIRED
\begin{abstract}
We introduce Hybrid-Cooperative Learning (HYCO), a framework for data-driven PDE modeling in which a physics-based solver and a flexible synthetic model are trained as two independent but cooperating components. Rather than imposing the governing equation as a residual on a single network, HYCO alternates between updating the physical parameters and the synthetic model, coupling them through agreement of their predictions.  Each component therefore fits the information available to it while acting as a regularizer for the other. This modular formulation accommodates sparse, heterogeneous, or disjoint datasets, avoids differentiating continuous PDE residuals, and combines standard solvers with general data-driven architectures. We further show that HYCO defines an exact potential game and establish equilibrium existence for a convex measure-relaxed model, providing a first structural interpretation of the alternating procedure. On inverse problems for reaction–diffusion systems, a heterogeneous Helmholtz equation, and a shock-forming traffic-flow model, HYCO reconstructs solutions and identifies parameters from sparse or localized observations, improving parameter recovery and extrapolation over uncoupled models, PINN-type methods, and classical solver-based inversion.
\end{abstract}


\maketitle

\section{Introduction}
\label{sec:intro}

\noindent
We introduce Hybrid-Cooperative Learning, HYCO in short, a hybrid learning framework designed to combine the strengths of physics-based and data-driven approaches while mitigating their respective limitations. In contrast to traditional methods such as PINNs \cite{RPK}, which impose physical and data-driven constraints directly on a single synthetic model, HYCO trains two models in parallel: one grounded in physical principles (the physical model), typically an ODE or PDE solver, and the other guided by data through a neural network ansatz (the synthetic model).

Our motivation stems from a fundamental modeling dilemma: the trade-off between the flexibility of purely data-driven models and the interpretability of physics-based models. This challenge has become even more pressing with the widespread adoption of machine learning (ML) techniques. Given a dataset, one typically faces a choice. A data-driven model is flexible and adaptive but unable to capture the underlying physics, whereas a physics-based model is rigid and interpretable, yet often difficult to specify and calibrate. Neither approach alone is sufficient to fully explain complex observations, especially when data is sparse or noisy. However, their strengths are naturally complementary: data-driven models excel in regions rich with observations, while physics-based models are better suited for extrapolation beyond the observed domain.

In HYCO, the two models interact and inform each other, much like two experts exchanging views before reaching a consensus on a complex problem. We formalize this cooperation through the following procedure:
\begin{itemize}
    \item[$\diamond$] We begin with a set of observations, whose availability may vary across models depending on the application. The two models may draw on data that is shared, partially overlapping, or entirely disjoint, including the limiting case where only the synthetic model is supplied with observations. The data may also be heterogeneous: one model might receive direct information about the solution, while the other relies on indirect measurements such as source terms, fluxes, or integral constraints. This flexibility in data allocation is a key feature of our approach.
    \item[$\diamond$] Training proceeds by alternating updates. At every step:
    \begin{itemize}
        \item The physical model receives the current prediction $\u_{syn}$ from the synthetic model and adjusts its parameters to fit the synthetic solution and the data it has access to (if any). Here, by ``parameters'' we mean the physical parameters appearing in the formulation of the PDE or ODE which is the physical model (e.g. Reynolds number, diffusivity, etc.).
        \item The synthetic model, based on a neural network ansatz, in turn, receives the prediction $\u_{phy}$ from the physical model and updates its parameters to match the physical prediction and the data it has access to. Here, the parameters are those appearing in the neural network formulation.
    \end{itemize}
\end{itemize}
This back-and-forth learning dynamic helps each model align not just with the data it sees and with the physical ansatz, but also with the behavior of its counterpart, leading to improved generalization, especially in regions where data is sparse.

\medskip
\noindent\textbf{Main contributions.}
The central contribution of this work is not a single model but a training framework. HYCO prescribes how a physics-based solver and a data-driven model can cooperate while remaining separate, rather than fixing what either of them is. The two are trained as distinct agents that communicate only through their predictions, each fitting whatever observations it can access while pulling the other toward agreement. Since the physical information reaches the synthetic model through the solver's output rather than a differential residual, HYCO never differentiates the governing equations and places almost no restriction on either component: the same scheme accepts any numerical solver in the physical slot and essentially any architecture in the synthetic one, which is precisely what makes it a reusable framework and not a one-off construction. We then put it to the test on three inverse problems of different character, elliptic (Helmholtz), reaction-diffusion (Gray-Scott), and hyperbolic (LWR), where the cooperative coupling recovers parameters and extrapolates more reliably from sparse or localized data than uncoupled models, PINN-type methods, or classical solver-based inversion. Finally, we read the framework through a game-theoretic lens: the alternating scheme is block-coordinate descent on the potential of an exact potential game, and in a convex relaxed setting that game admits a Nash equilibrium. This result is a structural justification, not a convergence guarantee for the algorithm used in the experiments, and we present it as such.

\subsection{Scope and Advantages of HYCO}
The scope of HYCO is broad. At its simplest, it can be used to embed standard machine learning models with physical biases, much like in the spirit of PINNs. The difference lies in how the physical bias is applied: instead of penalizing the PDE residual, the physical component is matched to the output of a numerical solver, whose discretization is locally accurate. Training against a computed solution in this way has been shown to be more accurate than training against the exact residual \cite{Esteve}.

More broadly, HYCO is naturally suited to merging partial information from heterogeneous sources. For instance, it can be applied in settings where we aim to reconcile observations in one region with a known physical model in another. A typical case is that of an incomplete PDE problem, where the right-hand side is known only in part of the domain, while solution measurements are available elsewhere. On this note, we remark that HYCO naturally applies to parameter identification and inverse problems, where mutual regularization improves recovery in settings that are ill-posed or lack sufficient information. This is a well-known challenge in the inverse problems literature and has become a major focus in scientific machine learning \cite{RCLV, WangWang}.

Compared to existing strategies, HYCO presents several further advantages:
\begin{enumerate}
    \item The method is apt to \textbf{distribute} the computational load. Indeed, the physical model and synthetic model can be implemented on different machines, thus splitting the computational cost of training. Besides, while the alternating scheme we adopt here is sequential, HYCO could theoretically be made fully \textbf{parallel} with minor changes (e.g. using a Jacobi variant of our scheme, in which both models are updated from the predictions of the previous iterate and swapped in together). This also opens the door to extensions with more than two models, which could have applications in domain decomposition or multiphysics contexts.
    \item Along the same lines, HYCO can be adapted to \textbf{privacy-aware} scenarios. Only the outputs of the models need to be exchanged, not their parameters or their structure. In fact, even the datasets could potentially differ, and different agents might collect data at different locations. As a result, each component can be trained locally with limited communication, mirroring the ideas behind federated learning \cite{FED}.
    \item Finally, when models can be shared, inference can be performed using either the physical or synthetic component, depending on the desired trade-off between accuracy and online speed of computation.
\end{enumerate}
Beyond these practical benefits, alternating descent has a structural advantage over the joint, single-loss training used by PINNs. In the PINN setting, a single network is fit simultaneously to the data and to the physical residual: both constraints act on the same set of parameters, and the joint gradient has no mechanism to distinguish a self-consistent but incorrect solution from the true one. Alternating updates break this symmetry: at each step, one model is optimized against a fixed prediction coming from an independent model family, rather than being free to reinterpret the physics to accommodate its own errors. We expect this to make it harder for the two models to drift into a shared, spurious optimum. We do not establish this mechanism theoretically, but it is consistent with the improved parameter recovery reported in Sections~\ref{sec:helm} and~\ref{sec:claw}, and it is particularly relevant for parameter identification problems, where the map from parameters to observations need not be uniquely invertible.

The most delicate point in HYCO concerns its computational cost. Indeed, the latter is inherently tied to the choice of the physical model: whenever the latter can be instantiated as a fast numerical solver, possibly even rooted in machine learning, the physical update is cheap, and training two models remains efficient overall. The bottleneck, in other words, is whether a numerical solver is available for the physical component, not the two-model structure per se.

\subsection{A Cooperative Modeling Game} HYCO can be easily framed in a game-theoretical perspective. This provides a theoretical justification to the approach, and helps in interpreting the inner workings of the methodology. In more technical terms, given two models with parameters $\Lambda$ (physical) and $\Theta$ (synthetic), we define three losses:
\begin{itemize}
    \item[$\diamond$] The physical loss $\L_{phy}(\Lambda)$, which is the distance between the physical model output, $\u_{phy}$, and the available observations. This term may also be omitted entirely if no data is assigned to the physical model.
    \item[$\diamond$] The synthetic loss $\L_{syn}(\Theta)$, which is the distance between the synthetic model output, $\u_{syn}$, and the observations;
    \item[$\diamond$] The interaction loss $\L_{int}(\Theta, \Lambda)$, which is a suitable distance between the prediction of the physical model and the prediction of the synthetic one.
\end{itemize}
Then, we aim to solve the following minimization problem
\begin{equation}
\label{main_minimization}
\begin{aligned}
&\min_{\Theta, \Lambda} &&\L := 
    \alpha\L_{syn}(\Theta) 
    +\beta \L_{phy}(\Lambda)
    + \L_{int}(\Theta, \Lambda) \\
&\text{subject to} \quad &&\u_{phy} \text{ solves the physical model} \\
&\, &&\u_{syn} \text{ solves the synthetic model},
\end{aligned}
\end{equation}
where the weights $\alpha, \beta >0$ serve as tunable hyperparameters to adjust the relative emphasis placed on each model.

\begin{figure}[t!]
  \centering
\includegraphics[width=.8\textwidth]{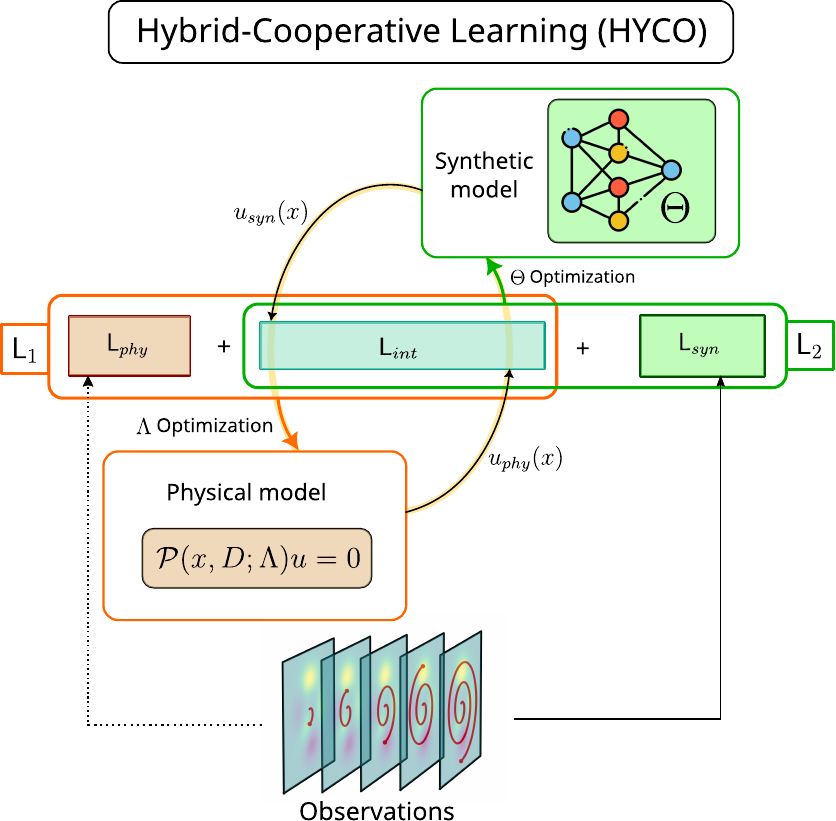}
  \caption{Schematic representation of the HYCO strategy. We train two models, alternating an optimization of the parameters $\Lambda$ of the physical model to one of the parameters $\Theta$ of the synthetic one. The goal of the physical model is to minimize its loss $\L_1$, which is the sum of $\L_{phy}$ and $\L_{int}$. Similarly, the synthetic model minimizes $\L_2$. The interaction loss $\L_{int}$ receives the two predictions $\u_{syn}, \u_{phy}$, while $\L_{phy}$ and $\L_{syn}$ are the discrepancies from the observations. The arrow leading from the observations to the physical loss is dotted because one might choose not to supply the physical model with data.}
  \label{fig:HYCO}
\end{figure}

In this context, by ``solves'' we mean that each model generates its output by adhering to its own internal structure: for the physical model, this typically involves solving an ODE or PDE; for the synthetic model, this may involve a static approximation, such as a neural network, a Gaussian mixture model, or even a least-squares fit, but also a learned dynamical system, such as a Neural ODE. A visual summary of the overall framework is provided in Figure~\ref{fig:HYCO}.

The best possible outcome is a Nash equilibrium: a pair of strategies $(\u_{phy}, \u_{syn})$ such that neither player can reduce their loss by unilaterally changing their parameters. In practice, especially when the true dynamics governing the data are unknown or inaccessible, this equilibrium embodies a compromise: by casting the problem as a game between these two imperfect models, we acknowledge their limitations and allow them to inform one another. The Nash equilibrium does not guarantee a perfect reconstruction of the data, but under uncertainty, it represents a stable and justifiable synthesis and, arguably, the best one can achieve given the circumstances. For further discussions on the matter, we refer to Section~\ref{sec:games}.

\subsection*{Outline of the paper}
After recalling the necessary notation, error metrics, and stopping criterion in Section~\ref{sec:error}, we begin with a concrete numerical experiment in Section~\ref{sec:grayscott}: a reaction-diffusion system exhibiting pattern formation. This serves as an illustrative case study to highlight the core concepts of HYCO. Building on insights from this example, we discuss the state of the art and the key advantages of our methodology in Section~\ref{sec:state}, and introduce the general abstract formulation, independent of any specific 
problem setup and applicable to a wide range of scenarios, in Section~\ref{sec:absstra}. We then turn to static settings, presenting an inverse problem for the Helmholtz equation in Section~\ref{sec:helm}, along with an in-depth sensitivity study, before returning to a dynamic setting with the LWR traffic model in Section~\ref{sec:claw}. Section~\ref{sec:games} explores the game-theoretic framework underpinning our approach and establishes the existence of a Nash equilibrium in a relaxed setting. We conclude with perspectives and future directions in Section~\ref{sec:conclusions}. A concrete recipe for instantiating this procedure on a new problem, together with practical guidance on the weights $\alpha, \beta$, the ghost-point count $H$, and the situations in which HYCO is most advantageous, is collected in Appendix~\ref{app:guide}.

\smallskip
\noindent
\textbf{Reproducibility.} The code reproducing all experiments is available at \href{https://https://github.com/ALive95/HYCO}{HYCO}. All experiments were run on a single MIG partition (1g.24gb, 24,GB) of an NVIDIA RTX PRO 6000 Blackwell Server Edition GPU. The implementation is based on JAX and Optax, with the finite element solves performed through a $P_1$ FEM discretization and Gmsh~\cite{gmsh} for mesh generation. Mesh sizes, time steps, network architectures, learning rates, ghost-point counts, and stopping-criterion parameters are reported in full within each experimental section.

\section{Mathematical Preliminaries}
\label{sec:error}

\noindent
Throughout this work we make repeated use of a few mathematical objects and evaluation metrics. This section collects their definitions in one place: the neural network architecture underlying our synthetic models, the error metrics used to compare HYCO against its competitors, and the stopping criterion adopted across all numerical experiments. Notationally, we use boldface for vector-valued quantities (e.g., $\x$, $\u$) and reserve non-bold symbols for scalars.

\subsection{Neural networks}
Deep neural networks (DNNs) are the synthetic models used throughout the numerical experiments of this paper. In most cases we equip them with residual connections, which we favor over classical feedforward DNNs since they ease optimization and improve gradient flow in deeper networks \cite{ResNet}; in some experiments, however, a plain feedforward architecture is sufficient and is used instead. A DNN is a parametric function $\u_{syn}$ that takes as input a vector $\x \in \R^p$ and returns a prediction of the ground truth at that point.

Formally, given a network with $\ell$ layers, $\u_{syn}$ is defined recursively as
\[
\u_{syn}(\x) = h_\ell,
\]
where
\begin{equation}
\label{eq:nn}
\begin{cases}
h_0 = \x \in \R^{p},\\
h_1 = \sigma(W_1 h_0 + b_1), \\
h_k = h_{k-1} + \sigma(W_k h_{k-1} + b_k) \text{ for } k=2,\ldots,\ell-1, \\
h_\ell = W_\ell h_{\ell-1} + b_\ell.
\end{cases}
\end{equation}
The residual connection is expressed by the presence of $h_{k-1}$ on the right-hand side of the third equation above; dropping this term recovers a standard feedforward network.
The function $\sigma$ is the activation function, which will be set to be either the ReLU function $\sigma(x) = \max\{x, 0\}$ or the hyperbolic tangent $\sigma(x) = \tanh(x)$, depending on whether smoothness is required or not.
The weight matrices and bias vectors are
\[
W_k \in \R^{q_k \times q_{k-1}}, \quad b_k \in \R^{q_k} \quad \text{for } k=1,\ldots,\ell,
\]
with $q_0 = p$ and $q_\ell = r$ denoting the input and output dimensions, respectively. We collect all trainable parameters in $\Theta = \{W_k, b_k\}_{k=1}^\ell$.

\subsection{Error metrics}
We compare HYCO against its competitors along three complementary metrics. Throughout, $m$ denotes an arbitrary model under consideration (e.g., HYCO's physical or synthetic component, a PINN, a plain DNN, etc.).

The most direct metric is the discrepancy between a model's prediction and the available data. Given a dataset $\D = \{\x_i, \u^\star(\x_i)\}_{i=1}^M$, where $\u^\star$ denotes the ground-truth function generating the data, we define
\begin{equation}
\label{eq:dataerr}
\mathsf{e}^m_d := \frac{1}{M}\sum_{i=1}^M \|\u_m(\x_i) - \u^\star(\x_i)\|^2 .
\end{equation}
This quantifies how well the model fits the observations it was trained on, but says nothing about its behavior elsewhere. To assess extrapolation, we compute the relative $L^2$ error between the predicted and true solution,
\begin{equation}
\label{eq:normL2}
\mathsf{e}^m_s := \frac{\|\u^\star-\u_m\|_{L^2}}{\|\u^\star\|_{L^2}}.
\end{equation}
Finally, for inverse problems in which hidden ground-truth physical parameters $\Lambda^\star$ must be recovered, we measure the relative error between the identified and true parameters,
\begin{equation}
\label{eq:normeu}
\mathsf{e}^m_p := \frac{\|\Lambda^\star - \Lambda_m\|}{\|\Lambda^\star\|}.
\end{equation}
Here $\|\cdot\|$ is the euclidean norm of the vector space of the physical parameters.

\subsection{Stopping criterion}
Comparing different models fairly requires accounting for their differing training requirements: each architecture converges after a different number of epochs, so a common, well-defined stopping rule is needed.

For parameter identification problems, we adopt an early stopping criterion based on the stabilization of the identified parameters, halting training once
\begin{equation}
\label{stopping}
\Big\|\Lambda^{(k)} - \frac{1}{Z} \sum_{j=k-Z}^{k-1} \Lambda^{(j)} \Big\| < \varepsilon,
\end{equation}
where $\Lambda^{(k)}$ denotes the identified parameters at epoch $k$, $Z$ is the averaging window size, and $\varepsilon$ is a fixed tolerance. The same rule is applied consistently across all numerical experiments in this paper, with $Z$ and $\varepsilon$ tuned per method.

\section{An Introductory Example: the Gray-Scott Model}
\label{sec:grayscott}

\noindent
As a testbed, we consider the well-known Gray-Scott reaction-diffusion system \cite{GrayScott}, which models autocatalytic and irreversible chemical reactions between two species, $u$ and $v$. This system is a prototypical example of nonlinear pattern formation, and its rich, stiff dynamics make it a well-known benchmark for assessing the capabilities of PDE learning models \cite{theWELL}. Purely data-driven models often struggle to capture the underlying dynamics, especially when asked to extrapolate beyond the observed data range. Physics-Informed Neural Networks (PINNs), which encode the governing equations as soft constraints within the loss function, tend to suffer from instability in this context, largely due to the stiffness of the equations and the scarcity of available data.

In this experiment, we evaluate HYCO's ability to learn the full solution, while identifying the model's diffusivity parameters at the same time. Consistent with the general framework introduced earlier, observational data is provided only to the synthetic model, while the physical model is trained solely through the interaction term (corresponding to taking $\beta = 0$ in $\L_1$).

On the square $\Omega = [0, 1] \times [0, 1]$, the system is defined as
\begin{equation}
    \label{gray_scott}
    \begin{dcases}
    u_t - D_u \Delta u + uv^2 - F(1 - u) = 0, \\
    v_t - D_v \Delta v - uv^2 + (F + k)v = 0, \\
    u(0,x,y) = u_0(x,y), \quad v(0,x,y) = v_0(x,y), \\
    \text{Periodic boundary conditions}.
    \end{dcases}
\end{equation}
In particular, $\u^\star = (u,v)$.
Here, $D_u$ and $D_v$ are the diffusion coefficients for $u$ and $v$, respectively; $F$ is the feed rate of $u$; and $k$ denotes the decay rate of $v$. The ground-truth parameters used to generate the data in this experiment are:
\[
D_u = 2\times 10^{-5}, \quad D_v = 8\times 10^{-6}, \quad F = 0.018, \quad k = 0.051.
\]
The initial values are
\[
u_0(x,y) = 1 - \tfrac12\,\chi_{B_{0.1}(\mathbf{c})}, \quad
v_0(x,y) = \tfrac12\,\chi_{B_{0.1}(\mathbf{c})},
\]
where $\chi_A$ denotes the characteristic function of the set $A$ and $B_{0.1}(\mathbf{c})$ is the ball of radius $0.1$ centered at $\mathbf{c} = (0.5,0.5)$.
In what follows, the initial datum is assumed known and is supplied to the physical model; the PINN competitor enforces it through an additional initial-datum penalty, while the purely data-driven network and the HYCO synthetic model receive no explicit initial-condition information beyond the dataset itself.

The dataset is constructed by numerically solving the Gray-Scott equations on a 64$\times$64 grid, using a finite difference scheme with small temporal steps ($5000$ time steps from $0$ to $T=2000$) to generate a reference true solution $\u^\star$. From this solution, we sample the values of both species at $M = 5000$ random space-time points, obtaining the dataset
\[
\D := \{\u^\star(x_i,y_i,t_i)\}_{i=1}^{M} = \{u^\star(x_i,y_i,t_i),v^\star(x_i,y_i,t_i)\}_{i=1}^{M}.
\]
Note, in particular, that our data is not equispaced in space nor in time.

The two models cooperating together in the HYCO strategy are:
\begin{itemize}
\item[$\diamond$] The physical model, governed by \eqref{gray_scott}, solved using a finite difference scheme with the same discretization used to generate the data. In this experiment, we supply the physical model of the initial datum $(u_0,v_0)$. However, the diffusivity coefficients $D_u, D_v$ are considered unknown. In the notation of Section~\ref{sec:error}, these are the parameters that the physical model will try to identify
\[
\Lambda^\star = \{ D_u, D_v\} = \{2\times 10^{-5},\, 8\times 10^{-6}\}.
\]
In HYCO and its competitors, these parameters are initialized far from the true values, at
\[
D_u = 1\times 10^{-5}, \quad \text{and} \quad D_v = 5\times 10^{-6}.
\]
\item[$\diamond$] The synthetic model is implemented as a feedforward deep neural network, a plain variant of \eqref{eq:nn} without the skip connections, with ReLU activation. The input is the vector $(x,y,t)$, and the network predicts the value of the solution at that point
\[\u_{syn}(x,y,t) = (u_{syn}(x,y,t), v_{syn}(x,y,t)).\]
In more detail, the neural network has four hidden layers, each containing 128 neurons.
We indicate the set of trainable parameters by $\Theta$. The synthetic loss is thus defined as
\[
    \L_{syn}(\Theta) := \frac1{M} \sum_{i=1}^M \|\u_{syn}(x_i,y_i,t_i) -  \u^\star(x_i,y_i,t_i)\|^2
    \]
\end{itemize}
Note that we do not penalize the parameters $\Theta$ of the synthetic model, as the current form of the functional, through its penalty on the distance to the physical approximation $\u_{phy}$, already ensures sufficient coercivity.

To complete the HYCO setting, we introduce an interaction loss that measures the discrepancy between the synthetic and physical models. To this end, at every epoch we select $H = 10^5$ random ``ghost points" $(x_i^G, y_i^G, t_i^G)$ (with the superscript ``G" standing for ``ghost"): the spatial locations are drawn uniformly in $\Omega$, while the time coordinates are drawn from the snapshot times stored by the physical solver, so that the physical prediction at each ghost point is an exact solver snapshot requiring no temporal interpolation. The ghost set is resampled afresh at every epoch. At these locations, we evaluate the predictions of both the physical and the synthetic models.
The interaction loss is then defined as:
\[
\L_{int}(\Theta, \Lambda) :=  \frac1{H} \sum_{i=1}^H \|\u_{syn}(x_i^G,y_i^G,t_i^G) - \u_{phy}(x_i^G,y_i^G,t_i^G)\|^2.
\]
This approach, very much in the spirit of Random Batch Methods (RBM), avoids the prohibitive cost of comparing the two predictions over the entire space-time grid at every iteration. In practice, as we shall see, evaluating the interaction loss at randomly resampled ghost points is sufficient to enforce physical consistency with significantly lower computational needs.

In summary, the synthetic model minimizes $\L_2 = \alpha\,\L_{syn} + \L_{int}$ over $\Theta$, for $\Lambda$ fixed, in an alternating fashion with the physical model's update. Here $\alpha > 0$ is a tunable hyperparameter balancing data fidelity and physical consistency in the synthetic loss.

To evaluate HYCO, we compare its performance to:
\begin{itemize}
    \item[$\diamond$] A purely data-driven feed-forward neural network, without the physical bias, implemented by setting the interaction loss in $\L_2$ to zero, effectively training only the synthetic model.
    \item[$\diamond$] A classical PINN, namely, a neural network function
    \[\u_{\text{PINN}}(x,y,t) = (\hat u(x,y,t),\hat v(x,y,t)).\]
    with parameters $\Theta$, optimized to minimize:
\[
\Tilde{\L}_{\text{PINN}}(\Theta,\Lambda) = \L_{\text{PINN}}(\Theta) + \gamma\mathsf{R}(\Theta, \Lambda) + \gamma_{0}\,\L_{0}(\Theta),
\]
where $\L_{\text{PINN}}(\Theta)$ is the discrepancy with respect to the data
\[
\L_{\text{PINN}}(\Theta) = \frac{1}{M}\sum_{i=1}^M \|\u_{\text{PINN}}(x_i,y_i,t_i) - \u^\star(x_i,y_i,t_i)\|^2,
\]
$\L_{0}(\Theta)$ penalizes the mismatch with the known initial datum on the grid at $t=0$, and $\mathsf{R}(\Theta,\Lambda)$ is the PDE residual
\begin{align*}
\mathsf{R}(\Theta,\Lambda) &= \|\hat u_{t} - D_u \Delta \hat u + \hat u \hat v^2 - F(1-\hat u)\|^2 \\
&\qquad + \|\hat v_t - D_v \Delta \hat v - \hat u \hat v^2 + (F + k)\hat v\|^2.
\end{align*}
The hyperparameters $\gamma, \gamma_0$ quantify the strength of each loss term.
Note that the PINN is also assumed not to know the diffusivity parameters $D_u, D_v$, hence the dependence of the residual term also on $\Lambda$.

In this experiment, the PINN consists of the same architecture of the HYCO neural network (four hidden layers with 128 neurons each), but uses the $\tanh$ activation. This is necessary because the PINN loss requires the activation function to be twice differentiable. The residual is evaluated at $10^4$ collocation points, drawn uniformly at random in the space-time domain at the start of training.
\end{itemize}

\smallskip
\noindent
We do not include a classical solver-based inverse method (e.g. adjoint or Gauss-Newton fitting of $D_u, D_v$) among the competitors here. By design, in this experiment the physical model receives no data ($\beta = 0$) and is driven solely by the interaction term; a solver-based inverse method, which fits the parameters directly against the observations, would therefore be solving a different problem, and the comparison would not be fair.

\smallskip
\noindent
All networks are trained with the Adam optimizer with learning rate $1 \times 10^{-3}$; the physical parameters are optimized in logarithmic scale with learning rate $8\times 10^{-3}$ in HYCO and $1\times 10^{-4}$ in the PINN. In HYCO, the synthetic network is first pretrained on the data alone for $2000$ epochs; the alternating phase is then run until the stopping criterion \eqref{stopping} is met (here with $Z = 1000$ and $\varepsilon = 10^{-4}$), which occurred after $5973$ epochs. The purely data-driven network, whose loss involves no identified parameters, is instead stopped once its loss stagnates, after roughly $8200$ epochs; the PINN, subject to the same rule \eqref{stopping} as HYCO, did not stop, and trained for the full 20000 epochs. In all cases the network is subsequently refined with up to $2000$ iterations of L-BFGS-B with the physical parameters frozen; for HYCO, this final refinement fits the synthetic network against the data and the physical prediction on the full grid at every stored snapshot. The strength parameter $\alpha$ in the HYCO loss is here set at $10$, while $\gamma$ and $\gamma_0$ in the PINN loss are both set at $100$. Since this section serves as an introductory illustration of the method, the results below are from a single run; systematic multi-seed statistics and wall-clock comparisons are carried out for the experiments of Sections~\ref{sec:helm} and~\ref{sec:claw}.

\begin{figure}[h!]
    \centering
    \includegraphics[width=0.8\linewidth]{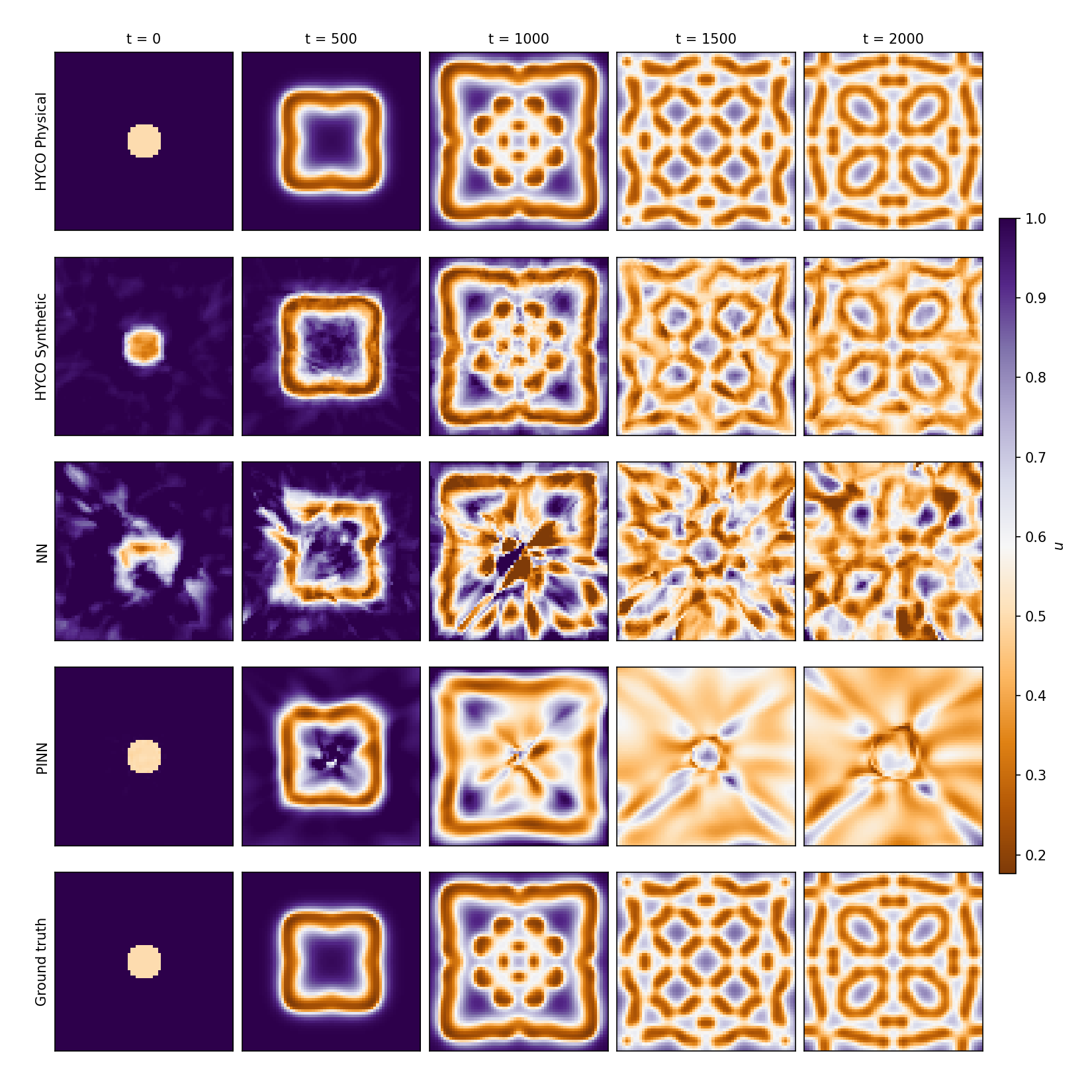}
    \caption{Comparison of model outputs for the Gray-Scott reaction-diffusion experiment, showing only the $u$-component of the solution. \textbf{First row:} Prediction from the HYCO physical model. \textbf{Second row:} Prediction from the HYCO synthetic model. \textbf{Third row:} Prediction from a purely data-driven neural network trained without physics-based guidance. \textbf{Fourth row:} Prediction from a Physics-Informed Neural Network (PINN). \textbf{Fifth row:} Ground-truth reference solution.}
    \label{fig:GS}
\end{figure}

\smallskip
\noindent
\textbf{Conclusions.}
\begin{itemize}
\item[$\diamond$] Figure~\ref{fig:GS} displays the evolution of the solution at five time snapshots. Every model captures the qualitative structure of the solution in the early phase of the dynamics. However, HYCO is the only method that remains accurate throughout the whole time horizon, losing little precision for the synthetic model only in the last snapshots. The purely data-driven network develops spurious small-scale oscillations as the pattern grows in complexity, while the PINN progressively over-smooths the solution for $t \gtrsim 1000$, as its limited physical structure is unable to keep up with the evolving pattern-formation dynamics.
\item[$\diamond$] In Table \ref{tab:gs}, we collect the final normalized errors between each model and the true solution, $\mathsf{e}^m_s$ computed as in \eqref{eq:normL2}.
The $L^2$ norm is computed on the same discretization grid used for the computation of the true solution.

We observe that HYCO Physical achieves the best precision, followed by HYCO Synthetic; both improve on the PINN and the purely data-driven network by a wide margin.
\item[$\diamond$] We also report, in Table \ref{tab:gsparams}, the identified coefficients by the PINN and HYCO Physical models. Clearly, we omit the HYCO synthetic model and the NN, as these models have no notion of physical parameters to reconstruct.
As we can see, HYCO is able to recover the ground truth coefficients, while the PINN struggles. A likely reason is that the solution network in the PINN is flexible enough to produce an acceptable approximation of the state even with biased parameters. Thus, during training, the optimizer tends to prioritize reducing the data mismatch, resulting in a solution field that looks reasonable, but with parameter values that deviate from the true one. Although we have experimented also with stronger residual enforcement (i.e. with larger $\gamma$ in the PINN loss), the solution network retains enough flexibility to absorb parameter errors, producing solutions that satisfy the PDE approximately while masking deviations in the coefficients.
\end{itemize}

\begin{table}[h!]
\centering
\renewcommand{\arraystretch}{1.2}
\begin{tabular}{|c|c|c|c|}
\hline
\textbf{HYCO Physical} & \textbf{HYCO Synthetic} & \textbf{PINN} & \textbf{NN} \\
\hline
\textbf{0.0146} & 0.0839 & 0.1832 & 0.2318 \\
\hline
\end{tabular}
\caption{Final normalized $L^2$ error for the approximations of the Gray-Scott equation.}
\label{tab:gs}
\end{table}

\begin{table}[h!]
\centering
\renewcommand{\arraystretch}{1.2}
\begin{tabular}{|c|c|c|}
\hline
& \textbf{$D_u$} & \textbf{$D_v$} \\
\hline
\textbf{Ground truth} & $2\times 10^{-5}$ & $8\times 10^{-6}$ \\
\hline
\textbf{HYCO Physical} & $\mathbf{2.004\times 10^{-5}}$ & $\mathbf{8.01\times 10^{-6}}$ \\
\hline
\textbf{PINN} & $1.151\times 10^{-5}$ & $5.84\times 10^{-6}$ \\
\hline
\end{tabular}
\caption{Identified diffusivity parameters $D_u, D_v$.}
\label{tab:gsparams}
\end{table}

\section{State of the Art}
\label{sec:state}

\noindent
Mathematical modeling traditionally relies on selecting a parametric structure, be it an algebraic equation, an ODE or PDE, or a neural network, and adjusting its parameters to best fit observed data. Two dominant paradigms have emerged: physics-based and data-driven modeling.

Physics-based models are grounded in physical laws and typically involve a small number of interpretable parameters. Their strength lies in theoretical guarantees, generalizability, and extrapolation. Classical examples include the Navier-Stokes equations or conservation laws in continuum mechanics, where the form of the equation is fixed and only a few parameters (e.g., Reynolds number, diffusion coefficients) need identification.

In contrast, data-driven models rely on flexible, high-capacity function approximators (typically neural networks) and large datasets. These models do not require prior knowledge of governing laws and have shown exceptional empirical performance, especially in unstructured tasks such as image and speech recognition \cite{MCGK}. However, their performance degrades in low-data regimes, and they tend to lack interpretability and extrapolation capacity.

\subsection{Hybrid Modeling}
To combine the complementary strengths of both paradigms, hybrid models have been introduced. The idea is simple: whenever some physical knowledge is available, one should exploit it to guide the learning process. This approach has gained substantial attention in recent years, especially in scientific machine learning, where data is expensive and often sparse \cite{BPK, RGQ}.

A well-known example is Physics-Informed Neural Networks (PINNs) \cite{RPK}, where the loss function penalizes not only the distance from the data but also from satisfying a known differential equation. For inverse and parameter-identification problems specifically, self-adaptive variants such as SA-PINN \cite{SAPINN} reweight the collocation loss to concentrate on high-residual regions, improving on vanilla PINNs at the cost of an additional adversarial weight update. Other approaches include Sparse Identification of Nonlinear Dynamics (SINDy), Neural Operators \cite{KLLABSA}, or Latent Dynamics Models \cite{RPSDQ}, which integrate physics either explicitly or in latent space representations. All these methods rely on one crucial assumption: an ansatz for the underlying physical law is available and can be enforced during training.

These strategies seek to find a compromise between model expressivity and structure, producing a single model that reflects both physical prior and data fit. However, enforcing hard constraints on the model may overly restrict its flexibility, especially when the physical knowledge is partial or noisy.

\subsection{The HYCO Perspective}
Rather than enforcing a single hybrid model, HYCO seeks to optimize two separate models, a physical and a synthetic one, penalizing their distance from data and their mutual prediction discrepancies.  

This architecture draws inspiration from ideas already present in the literature. A first connection is to Ensemble Learning \cite{BO, ZM}, where predictions from multiple models are aggregated together. However, the key difference is that our aggregation is implicit and done during training, not after. In a similar vein, one can view the interaction term as a regularizer, analogous to PDE-based regularization \cite{AANS, NM}, where the synthetic model is biased toward physically consistent behavior by the proximity to the physical model. 
The HYCO method also closely relates to recent ideas used in ADMM-PINN hybrids \cite{YXH}, where different models are coordinated through consistency terms. Yet, our method makes no assumption of model hierarchy or structure: the two models are independent in architecture and are trained solely through consistency with data and with each other. Unlike PINNs or ADMM-PINNs, however, HYCO does not modify the synthetic architecture to enforce the physics, nor makes assumptions on the model hierarchy. Indeed, both models influence each other iteratively and indirectly.
Finally, the game-theoretical framework behind HYCO shares conceptual similarities with the Federated Learning paradigm \cite{FED}, where multiple agents collaborate to train a model without exchanging data. Although HYCO does not follow the federated approach in a strict sense, since it involves the cooperative training of distinct models rather than a shared one, it retains a similar philosophy. Each agent operates independently and communicates only through model predictions, preserving the privacy. Connections between these distributed learning settings and game theory are increasingly explored in the literature; see for instance \cite{MMRHA, WaZu}.

Two further lines of work couple a solver and a network in a related but distinct manner. WF-PINNs \cite{WFPINN} pair a main and an auxiliary network through a consistency constraint for inverse problems, but both remain soft, PINN-style approximators rather than one of them being an explicit numerical solver. PRDP \cite{PRDP} instead embeds an iterative differentiable solver directly in the computational graph, adaptively refining its accuracy during training to reduce cost; unlike HYCO, solver and network are optimized jointly through a single objective rather than through alternating, mutually-regularizing updates.

\section{The Abstract HYCO Strategy}
\label{sec:absstra}

\noindent
In this section we introduce the abstract HYCO strategy, formulated in continuous terms to retain full generality and remain independent of any specific discretization scheme. This formulation ensures that the framework can be flexibly adapted to a range of numerical methods. For example, the discretizations adopted in Sections~\ref{sec:grayscott} and~\ref{sec:helm} are only one possible choice; others may lead to different results.

The HYCO methodology is fully agnostic to the specific model structure: given access to a physical ansatz and a synthetic (data-driven) model, HYCO integrates them into a unified training framework. The core ingredients are:
\begin{itemize}
    \item[$\diamond$] Observational data sampled from the ground-truth system;
    \item[$\diamond$] A physical ansatz, namely a PDE or ODE system assumed to govern the phenomenon, depending on unknown structural parameters $\Lambda$;
    \item[$\diamond$] A purely data-driven synthetic model depending on parameters $\Theta$.
\end{itemize}
Once these components are defined, the objective is to jointly optimize $\Lambda$ and $\Theta$ to fit the data. Since $\Lambda$ encodes physical parameters, this automatically yields a solution to the corresponding inverse problem. We stress that the physical ansatz and the physical model are distinct: the former is the governing equation itself, while the latter is whatever means HYCO uses in practice to approximate a solution of it. This is typically a classical numerical discretization (finite differences, finite elements, finite volumes, spectral methods), but it may in principle also be a neural network or operator, provided it is constructed to approximate the ansatz rather than being fit directly to data as the synthetic model is.

\smallskip
\noindent
\textbf{A word of warning.}
We assume, for simplicity, that both models have access to the full dataset. As noted in the Introduction and the Gray-Scott experiment, this assumption is not required in practice, since the models only need to share their predictions.

We present the abstract strategy in the dynamic setting, as applied in Section~\ref{sec:grayscott}  and~\ref{sec:claw}. The static setting discussed in Section~\ref{sec:helm} then emerges as a special case of this general framework.

\smallskip
\noindent
\textbf{The Data.} To maintain maximum generality, we assume access to $M$ space-time measurements, collected at points $(\x_i, t_i) \in \R^d \times [0,T]$, with $T > 0$. The dataset is thus given by
\[
\mathcal{D} := \{\u^\star(\x_i, t_i)\}_{i=1}^{M},
\]
where $\u^\star(\x_i, t_i) \in \R^k$ is the measurement recorded at location $\x_i$ and time $t_i$.

This setting encompasses several sampling strategies. In the \emph{static} case, the data are collected by $M'$ fixed sensors at $N$ common time instants $\{t_j\}_{j=1}^N$, so that $\mathcal{D}$ takes the tensor-product form $\{\u^\star(\x_i, t_j)\}_{i=1, j=1}^{M', N}$ with $M = M'N$; this configuration is commonly encountered in both theoretical and applied contexts. In the \emph{mobile} case, the sensor locations vary in time, $\x_i = \x_i(t_i)$. More generally, the sampling points may be scattered in space-time with no common time grid at all, as in the Gray-Scott experiment of Section~\ref{sec:grayscott} and the LWR experiment of Section~\ref{sec:claw}, where the points $(\x_i, t_i)$ are drawn at random. All the constructions below use only the points $(\x_i, t_i)$ themselves, and therefore apply verbatim in each of these cases.

\smallskip
\noindent
\textbf{The Physical Model.} The physical ansatz is the PDE system
\begin{equation}
\label{dynamicPDE}
\begin{dcases}
\partial_t \u_{phy}(\x,t) + \mathscr{F}(\kappa)\u_{phy}(\x,t) = f, & \x \in \Omega,\ t \in (0,T], \\
\mathscr{G}(\lambda)\u_{phy}(\x,t) = 0, & \x \in \partial\Omega,\ t \in (0,T], \\
\u_{phy}(0,\x) = \u_0(\x), & \x \in \Omega,
\end{dcases}
\end{equation}
which governs the evolution of the unknown field $\u_{phy} : [0,T] \times \Omega \to \R^k$. The operators $\mathscr{F}(\kappa)$ and $\mathscr{G}(\lambda)$ are differential operators depending on structural parameters $\kappa$ and $\lambda$, which may be constants or functions. The source term $f$ and the initial condition $\u_0$ may be known or unknown. All unknown quantities are grouped into the set of physical parameters:
\[
\Lambda := \{\kappa(\x), \lambda(\x), \ldots\}.
\]
To obtain $\u_{phy}$ in practice, HYCO approximates this ansatz with a physical model: a numerical discretization, as in Sections~\ref{sec:grayscott}, \ref{sec:helm}, and \ref{sec:claw}, or, more generally, any solver, including a neural network or operator, provided it is constructed to satisfy the ansatz rather than being fit to data as the synthetic model is. If the physical model has access to data, it seeks to minimize the loss
\begin{equation}
\label{loss_phy}
\LL_{phy}(\Lambda) := \frac{1}{M} \sum_{i=1}^M \ell_{phy}(\u_{phy}(\x_i, t_i),\u^\star(\x_i, t_i))  + \mathscr{P}(\Lambda),
\end{equation}
subject to $\u_{phy}$ solving \eqref{dynamicPDE}. The function $\ell_{phy}$ is a positive bivariate function that quantifies the mismatch between the model and the data, typically the mean squared error in our experiments. The term $\mathscr{P}(\Lambda)$ is a regularization on the physical parameters. Here and below, $\LL$ denotes the abstract, continuous-in-time formulation of each loss; Sections~\ref{sec:grayscott}, \ref{sec:helm}, and \ref{sec:claw} recover the discrete, practical version $\L$ used in each algorithm, once a discretization or a ghost-point sampling scheme is fixed.

\smallskip
\noindent
\textbf{The Synthetic Model.} This component can take the form of any machine learning architecture. In this work, we consider the case of a feed-forward neural network $\u_{syn}$, as in \eqref{eq:nn}. The loss mirrors that of the physical model:
    \begin{equation}
    \label{loss_syn}
    \LL_{syn}(\Theta) := \frac{1}{M} \sum_{i=1}^M  \ell_{syn}(\u_{syn}(\x_i, t_i),\u^\star(\x_i, t_i)) + \mathscr{P}(\Theta),
    \end{equation}
    where $\mathscr{P}(\Theta)$ is a regularization term, typically used to promote properties such as sparsity.

\bigskip
\noindent
\textbf{The Interaction Loss.} 
The interaction loss measures the discrepancy between the synthetic and physical predictions over space-time. Accordingly, we simply set
\begin{equation}
\label{loss_int}
\LL_{int}(\Theta, \Lambda) := \int_0^T\int_\Omega \ell_{int}(\u_{syn}(\x,t),\u_{phy}(\x,t))\,d\x\,dt.
\end{equation}
Approximating numerically this term precisely at every epoch can make the training expensive. Therefore, in concrete examples, we decided to adopt the idea of ghost points, already introduced in Section~\ref{sec:grayscott}. Drawing from the theory of Random Batch Methods (RBM), at every epoch we sample a certain number $H$ of ghost points and compute the norm of the difference of the two models on this chosen set.

The overall procedure is summarized in Algorithm~\ref{alg:hyco}. The two updates are performed in Gauss-Seidel fashion: the synthetic update uses the physical parameters $\Lambda^{k+1}$ just produced within the same step. We stress that each block update is a single gradient step of the chosen optimizer, not an exact minimization of the corresponding loss; this distinction matters when interpreting the scheme through the game-theoretic lens of Section~\ref{sec:games}, and we return to it there.

\begin{algorithm}[h!]
\caption{HYCO alternating training}
\label{alg:hyco}
\KwIn{data $\D$, physical solver, synthetic model, weights $\alpha, \beta$, ghost count $H$}
\KwOut{$\u_{phy}, \u_{syn}, \Lambda$}
Initialize $\Lambda^0, \Theta^0$\;
\For{$k = 0, 1, 2, \dots$}{
  Sample $H$ ghost points\;
  $\Lambda^{k+1} \gets$ optimization step on $\L_1(\Theta^k, \Lambda) = \beta\,\L_{phy}(\Lambda) + \L_{int}(\Theta^k, \Lambda)$ with respect to $\Lambda$\;
  $\Theta^{k+1} \gets$ optimization step on $\L_2(\Theta, \Lambda^{k+1}) = \alpha\,\L_{syn}(\Theta) + \L_{int}(\Theta, \Lambda^{k+1})$ with respect to $\Theta$\;
  \If{stabilization criterion \eqref{stopping} met}{\textbf{break}\;}
}
\end{algorithm}

\section{A Static Example: The Helmholtz Equation}
\label{sec:helm}

\noindent
While we began by illustrating our algorithm on a dynamic problem, HYCO is equally suited to static settings. To demonstrate this, we consider the heterogeneous Helmholtz equation, a classical benchmark for both forward and inverse problems in the machine learning and PDE literature \cite{BPINNs, Fang, osti}. We pose the problem on a \emph{nonconvex} domain $\Omega \subset \R^2$, whose irregular geometry is naturally handled by the finite element discretization.

\subsection{Problem setup and HYCO components}

The ground-truth field $u^\star \in H^1_0(\Omega)$ solves the heterogeneous Helmholtz problem with homogeneous Dirichlet boundary conditions
\begin{equation}
\label{helmholtz}
\begin{dcases}
- \nabla \cdot (\kappa(\x) \nabla u^\star) + \eta(\x)^2\, u^\star = f, & \x \in \Omega, \\
u^\star = 0, & \x \in \partial\Omega,
\end{dcases}
\end{equation}
with forcing $f(\x) = 10\sin(x_1)\cos(x_2)$. The diffusion coefficient $\kappa$ and the wave number $\eta$ are modeled as unit-shifted Gaussian bumps,
\begin{equation}
\label{kappaeta}
\kappa(\x;\alpha_1,c_1) = 1 + \alpha_1 e^{-|\x - c_1|^2}, \qquad
\eta(\x;\alpha_2,c_2) = 1 + \alpha_2 e^{-|\x - c_2|^2},
\end{equation}
with ground-truth parameters collected in
\[
\Lambda^\star := \{\alpha_1, c_1, \alpha_2, c_2\} \in \R^6.
\]
Since $c_1, c_2$ are 2D centers, the total number of parameters to identify is six. The reference solution $u^\star$ has no closed form and is obtained by a high-fidelity $P_1$ FEM solve: the domain boundary is discretized with $80$ spline control points and meshed with Gmsh \cite{gmsh} at target element size $h_{fine} = 0.05$, giving a fine triangular mesh that accurately resolves both the geometry and the Gaussian bumps in $\kappa, \eta$.

We generate the data by placing $M = 25$ sensors at scattered interior locations and evaluating the solution there:
\[
\D = \{u^\star(\x_i)\}_{i=1}^{M}.
\]
The basic and comparison experiments below use clean observations; robustness to observational noise is examined separately in the sensitivity study of Section~\ref{sec:helm_sens}. 

\smallskip
\textbf{Goal:} {\it Given the dataset $\D$, the source $f$, and the model form \eqref{helmholtz}, identify the parameters $(\alpha_j, c_j)$ for $j = 1, 2$ such that the corresponding solution best matches the observations.}

\smallskip
\noindent
As in the general framework, HYCO couples two models.
\begin{itemize}
\item[$\diamond$] The \textbf{physical model} is the $P_1$ FEM discretization of \eqref{helmholtz} on a \emph{coarse} triangular mesh. In more detail, we set the target element size to $h_{coarse} = 0.2$, four times coarser than the mesh used to generate the data. The output of the physical model is the nodal solution $u_{phy}$ for a given parameter set $\Lambda$. We supply it with the data, so that its loss is
\[
\L_{phy}(\Lambda) = \frac{1}{M}\sum_{i=1}^{M} \|u_{phy}(\x_i) - u^\star(\x_i)\|^2 .
\]
\item[$\diamond$] The \textbf{synthetic model} is a feed-forward neural network $u_{syn}$, taking $\x = (x_1,x_2)$ as input and predicting the solution value, with $2$ hidden layers of width $256$, ReLU activation, and residual connections, as in \eqref{eq:nn}. Its loss is
\[
\L_{syn}(\Theta) = \frac{1}{M}\sum_{i=1}^{M} \|u_{syn}(\x_i) - u^\star(\x_i)\|^2 .
\]
\item[$\diamond$] The \textbf{interaction term} measures the discrepancy between the two models away from the data. We draw $H$ ghost sensors $\{\x_h^G\}$ in $\Omega$ at each epoch and set
\[
\L_{int}(\Theta,\Lambda) = \frac{1}{H}\sum_{h=1}^H \|u_{syn}(\x_h^G) - u_{phy}(\x_h^G)\|^2 .
\]
\end{itemize}
The physical model minimizes $\L_1 = \beta\L_{phy} + \L_{int}$ over $\Lambda$, and the synthetic model minimizes $\L_2 = \alpha\L_{syn} + \L_{int}$ over $\Theta$, in the alternating fashion described in Section~\ref{sec:absstra}. We adopt the stabilization-based stopping criterion \eqref{stopping} for the parameter-identifying methods.

\smallskip
\noindent
\textbf{Numerical parameters.} All models are trained with the Adam optimizer. In HYCO the physical and synthetic models use learning rates $5\times10^{-3}$ and $1\times10^{-3}$, respectively. The number of ghost points is $H = 200$ throughout training. We set $\alpha=10$ and $\beta=1$, so both models receive a direct data term in addition to the interaction loss. The PINN and SA-PINN baselines use the same network width, $\tanh$ activation (required for the residual autodiff), and $2000$ collocation points. Training stops once the stabilization criterion \eqref{stopping} is met, with window $Z=1000$ and tolerance $\varepsilon=10^{-3}$ for HYCO, PINN, and SA-PINN, and $Z=5$, $\varepsilon=10^{-4}$ for Gauss-Newton. All results reported in this section are averaged over five random seeds $\{42,33,17,23,2\}$ and given as mean $\pm$ standard deviation.

\subsection{Basic experiment}

\begin{figure}[h!]
    \centering
    % PLACEHOLDER
    \includegraphics[width=0.7\linewidth]{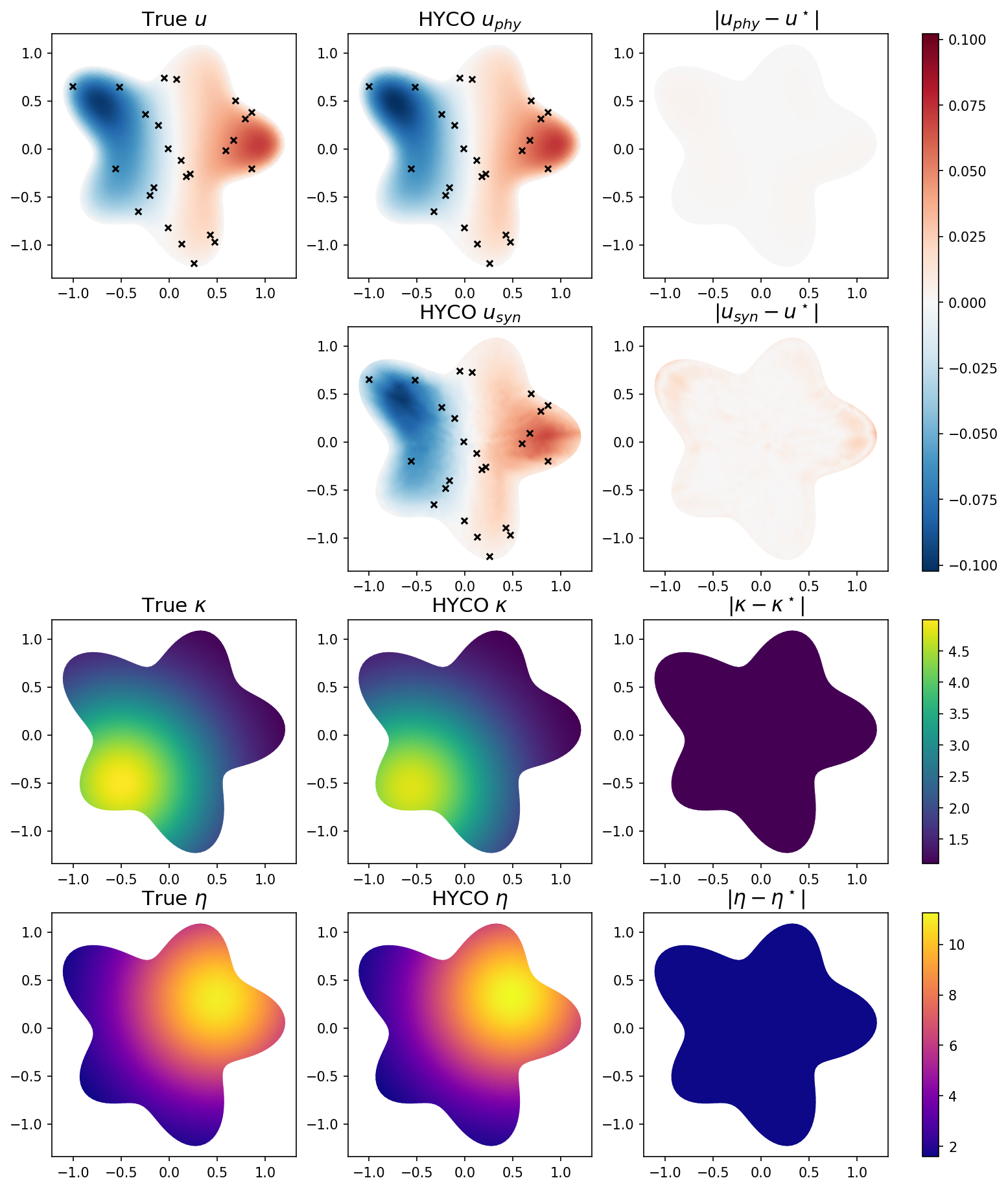}
    \caption{Helmholtz basic experiment. \textbf{Row 1:} reference solution $u^\star$, HYCO's physical prediction $u_{phy}$, HYCO's synthetic prediction $u_{syn}$, and the pointwise error $|u_{phy} - u^\star|$; black $\times$ markers indicate the sensor locations used for training. \textbf{Row 2:} ground-truth coefficient fields $\kappa^\star,\eta^\star$ and the corresponding HYCO reconstructions $\hat\kappa,\hat\eta$. Fields shown for seed $42$.}
    \label{fig:helm_basic}
\end{figure}

We first run HYCO in isolation to verify that it solves the inverse problem. Figure~\ref{fig:helm_basic} reports the reference solution, the recovered physical and synthetic fields, and the recovered coefficients $\hat\kappa$, $\hat\eta$ against the ground truth. Despite the sparse observations, both HYCO predictors reproduce $u^\star$ across the nonconvex domain, and the recovered coefficient fields match the ground-truth Gaussian bumps in location and amplitude.

\subsection{Ablation and comparison}

\noindent
We now turn to comparing HYCO against alternative approaches to the same inverse problem. This addresses two separate questions: how much of the method's performance comes from the cooperative coupling itself, and how HYCO stacks up against existing inversion strategies.

\smallskip
\noindent
\textbf{Ablation.} We isolate the contribution of the cooperative coupling by removing $\L_{int}$ from both losses. Without it, the synthetic network is trained on the sensor data alone, and the FEM solver reduces to a classical parameter-fitting inverse problem: the two models never exchange information at the ghost points, and each solves its own inverse problem independently. Table~\ref{tab:helm_comp} (rows HYCO-abl and HYCO) reports the resulting parameter error. The coupled HYCO improves parameter recovery over the decoupled baseline, confirming that the synthetic model provides a useful regularization signal to the physical inversion.

\begin{figure}[h!]
    \centering
    % PLACEHOLDER
    \includegraphics[width=\linewidth]{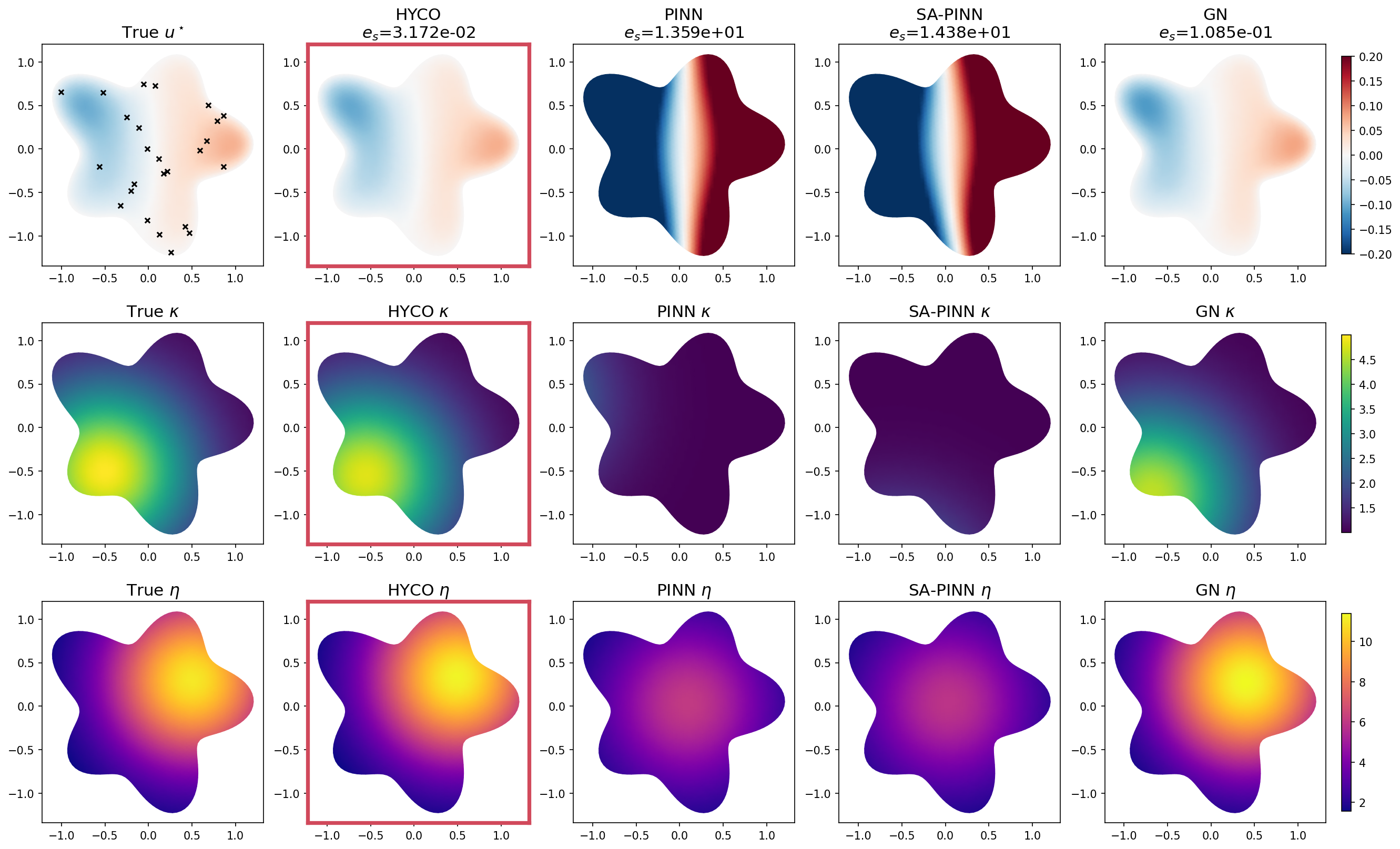}
    \caption{Helmholtz comparison. \textbf{Row 1:} the reference solution $u^\star$ (black $\times$ marking the sensor locations), followed by the recovered solution field for each competing method, with the corresponding solution error $\mathsf{e}^m_s$ \eqref{eq:normL2} reported above each panel; ``HYCO'' denotes the physical model's prediction $u_{phy}$. \textbf{Rows 2--3:} the ground-truth coefficient fields $\kappa^\star$ and $\eta^\star$, followed by the fields recovered by each method. Fields shown for seed $42$.}
    \label{fig:helm_comp}
\end{figure}

\smallskip
\noindent
\textbf{Competitors.} We compare HYCO against three benchmarks, chosen to span the range of existing approaches to this inverse problem: a standard physics-informed baseline, a stronger variant of it built specifically for inverse problems, and a purely classical solver-based method. We exclude Neural Operator baselines (e.g., FNO, DeepONet): these methods are trained across many instances of a parametric PDE family to learn a solution operator, whereas this experiment is a single-instance inverse problem with one fixed, unknown parameter set to recover, for which they are not designed.
\begin{enumerate}[label=(\roman*)]
\item \textbf{PINN} \cite{RPK}: a $\tanh$ network trained with $\L_{data} + \gamma\, \L_{pde}$, where $\L_{pde}$ penalizes the strong-form residual of \eqref{helmholtz} via automatic differentiation, with $\Lambda$ learned jointly. This is the default neural baseline for PDE-constrained inverse problems.
\item \textbf{SA-PINN} \cite{SAPINN}: a self-adaptive PINN maintaining a learnable weight vector over the collocation pool, updated by gradient ascent to concentrate on high-residual regions. Vanilla PINNs are known to underperform on inverse problems, and self-adaptive collocation weighting is the most established remedy, making this a stronger neural baseline than (i).
\item \textbf{Gauss-Newton (GN)}: a classical solver-based inversion using the update $\delta\Lambda = -(J^\top J + \lambda I)^{-1} J^\top r$, where $J = \partial G/\partial\Lambda \in \R^{M\times 6}$ is the Jacobian of the FEM forward map $G:\Lambda \mapsto \u_\Lambda(\x_{1:M})$ and $r = G(\Lambda) - u^\star$. This isolates how much of HYCO's advantage comes from having a learned surrogate at all, as opposed to any property specific to the FEM solver.
\end{enumerate}
Results are collected in Table~\ref{tab:helm_comp}. For the models trying to recover also the hidden coefficients, we show the recovered fields in Figure~\ref{fig:helm_comp}.

\begin{table}[h!]
\centering
\small
\renewcommand{\arraystretch}{1.2}
\setlength{\tabcolsep}{4pt}
\begin{tabular}{|c|c|c|c|c|}
\hline
\textbf{Method} & $\mathsf{e}^m_p$ & $\mathsf{e}^m_s$ & $\mathsf{e}^m_d$ & \textbf{ms/epoch\,$\vert$\,iter} \\
\Xhline{1.2pt}
PINN & \makecell{0.5132 \\ {\scriptsize $\pm0.0078$}} & \makecell{13.80 \\ {\scriptsize $\pm0.51$}} & \makecell{1.7971e-01 \\ {\scriptsize $\pm1.92$e-02}} & 52.8 \\
\hline
SA-PINN & \makecell{0.5335 \\ {\scriptsize $\pm0.0021$}} & \makecell{14.89 \\ {\scriptsize $\pm0.35$}} & \makecell{2.1224e-01 \\ {\scriptsize $\pm3.02$e-02}} & 51.6 \\
\hline
Gauss-Newton & \makecell{0.0611 \\ {\scriptsize $\pm0.0166$}} & \makecell{0.0893 \\ {\scriptsize $\pm0.0155$}} & \makecell{9.9714e-06 \\ {\scriptsize $\pm3.07$e-06}} & 123.4$^\dagger$ \\
\hline
HYCO-abl (phy, $\L_{int}=0$) & \makecell{0.0611 \\ {\scriptsize $\pm0.0166$}} & \makecell{0.0893 \\ {\scriptsize $\pm0.0155$}} & \makecell{9.9724e-06 \\ {\scriptsize $\pm3.07$e-06}} & 7.6 \\
HYCO-abl (syn, $\L_{int}=0$) & -- & \makecell{2.2030 \\ {\scriptsize $\pm0.5230$}} & \makecell{8.4660e-08 \\ {\scriptsize $\pm1.44$e-07}} & -- \\
\hline
HYCO (phy, ours) & \makecell{\textbf{0.0550} \\ {\scriptsize $\pm0.0365$}} & \makecell{\textbf{0.0546} \\ {\scriptsize $\pm0.0203$}} & \makecell{3.7573e-06 \\ {\scriptsize $\pm2.19$e-06}} & 9.6 \\
HYCO (syn, ours) & -- & \makecell{0.1474 \\ {\scriptsize $\pm0.0162$}} & \makecell{\textbf{1.2743e-08} \\ {\scriptsize $\pm4.28$e-09}} & -- \\
\hline
\end{tabular}
\caption{Helmholtz: ablation and comparison. Parameter error $\mathsf{e}^m_p$ \eqref{eq:normeu}, solution error $\mathsf{e}^m_s$ \eqref{eq:normL2}, and data discrepancy $\mathsf{e}^m_d$ \eqref{eq:dataerr} at the sensor locations. Best accuracy in bold. Error entries are mean $\pm$ standard deviation over five random seeds $\{42,33,17,23,2\}$. The last column is the per-step wall-clock cost on common hardware: milliseconds per epoch for the gradient-trained models, and per iteration for Gauss-Newton; the two HYCO / HYCO-abl predictors come from one coupled run. The HYCO-abl (phy) and Gauss-Newton rows coincide up to optimizer tolerance: the decoupled physical fit minimizes the same FEM misfit that Gauss-Newton solves, and both converge to the same minimizer on every seed. $^\dagger$Gauss-Newton is not epoch-based and converges in only $50$ iterations ($6.2$\,s total), so its per-iteration cost is not directly comparable to an epoch.}
\label{tab:helm_comp}
\end{table}

\noindent
\textbf{Conclusions.}
\begin{itemize}
\item[$\diamond$] All methods drive the data discrepancy $\mathsf{e}^m_d$ to small values, but they converge to different points in parameter space. HYCO recovers $\Lambda^\star$ most accurately, while PINN and, to a lesser extent, SA-PINN settle on alternative parameter sets that fit the data but deviate from the true coefficients. Interestingly, while the physical model achieves the best solution error, it is actually the synthetic model that attains the smallest data discrepancy. This is expected: unconstrained by the physics, the synthetic network can fit the sparse sensor values almost exactly, whereas the physical prediction must remain a valid FEM solution and therefore trades a slightly larger data discrepancy for markedly better accuracy away from the sensors. 
This complementarity is also practical: at inference one is free to query whichever model suits the task, the synthetic surrogate for cheap pointwise evaluation, or the physical model for a solution that stays consistent with the governing equations away from the data.
\item[$\diamond$] Gauss-Newton relies on a local linearization of the forward map and has no mechanism to regularize against the nonconvexity of the misfit landscape; HYCO's alternately co-trained synthetic model supplies exactly this kind of implicit regularization.
\item[$\diamond$] On the nonconvex geometry, HYCO does not require a fine physical mesh to be effective: each HYCO epoch costs about $9.6$\,ms, roughly five times less than the PINN and SA-PINN baselines ($\sim\!52$\,ms), whose per-epoch cost is dominated by differentiating the strong-form residual through the network. The cooperative coupling adds only about $30\%$ over the decoupled single-model fit ($7.6$\,ms per epoch), so the mutual regularization is nearly free per step. Gauss-Newton is the cheapest method overall in absolute terms, converging in $50$ iterations and about $6$\,s; its per-iteration Jacobian assembly ($123$\,ms) is heavier than an epoch, but the iteration count is small. It nonetheless settles on a less accurate parameter set than HYCO, and this is not a matter of budget: the misfit stalls after roughly $20$ iterations (Figure~\ref{fig:helm_comp}), so allowing more iterations does not close the gap. The limitation is the nonconvexity of the misfit landscape, against which HYCO's co-trained surrogate provides a regularizing signal that a purely local linearization lacks.
\end{itemize}

\subsection{Sensitivity study}
\label{sec:helm_sens}

We close the Helmholtz experiment by probing how HYCO depends on its three main
hyperparameters: the loss weights $\alpha,\beta$, the number of ghost points $H$,
and the level of observational noise. All studies use HYCO alone, since these
knobs have no analogue in the baselines. The non-varied weights are held at their
defaults $\alpha=10$, $\beta=1$ and $H=200$.

\begin{figure}[h!]
    \centering
    \includegraphics[width=\linewidth]{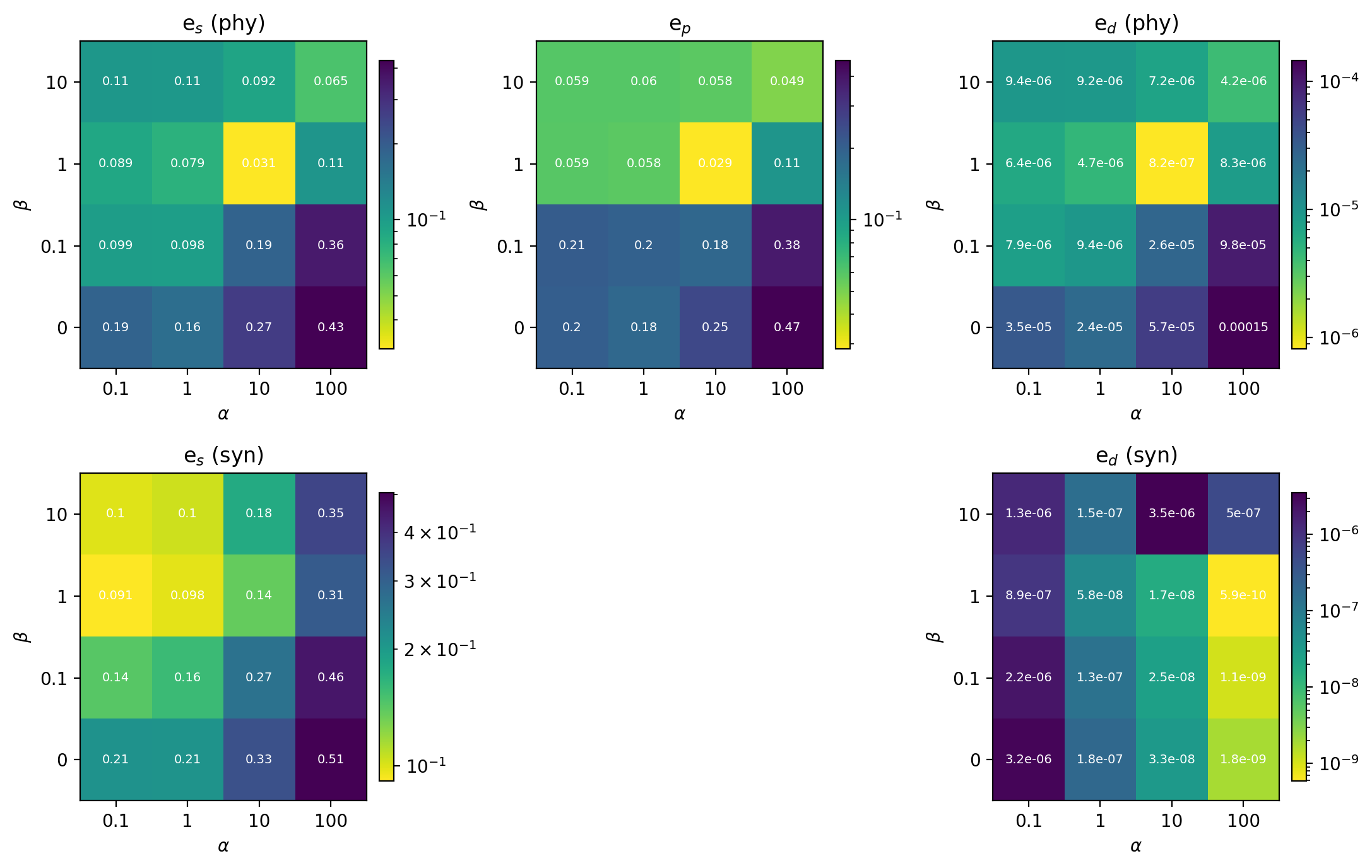}
    \caption{Helmholtz sensitivity to the loss weights. Error matrices over the
    $(\alpha,\beta)$ grid: physical solution error $\mathsf{e}^m_s$ and data
    discrepancy $\mathsf{e}^m_d$ and parameter error $\mathsf{e}^m_p$ (top row),
    and the synthetic-model solution error and data discrepancy (bottom row). The
    default $\alpha=10$, $\beta=1$ sits in the low-error interior plateau.}
    \label{fig:helm_ab}
\end{figure}

\smallskip\noindent
\textbf{Loss weights.} We sweep 
\[\alpha \in \{0.1, 1, 10, 100\}, \qquad \beta \in \{0, 0.1, 1, 10\}.\]
Recall that $\beta$ scales the physical data term and
$\alpha$ the synthetic one, while the interaction loss carries unit weight in both
objectives. Figure~\ref{fig:helm_ab} summarizes the errors into a $4\times4$ grid matrix. 

When $\beta=0$ the physical model
receives no data and is driven by the interaction term alone. Unsurprisingly, in this setting parameter recovery degrades by roughly an order of magnitude ($\mathsf{e}^m_p$ between $0.18$ and $0.47$). At the opposite extreme, $\alpha=100$ over-weights the
synthetic data fit: the synthetic discrepancy $\mathsf{e}^m_d$ falls to
$\sim\!10^{-9}$ while its solution error $\mathsf{e}^m_s$ rises to $0.31$--$0.51$,
the signature of a network that interpolates the $25$ sensors at the expense of
everywhere else. Between these extremes the physical metrics are flat: over
$\alpha\in[1,10]$, $\beta\in[0.1,1]$ the parameter error stays near its minimum,
attained at the default $\alpha=10$, $\beta=1$ ($\mathsf{e}^m_p=0.029$,
$\mathsf{e}^m_s=0.031$). Performance is therefore not knife-edge in the weights:
any choice in this interior plateau gives comparable accuracy.

\begin{figure}[h!]
    \centering
    \includegraphics[width=0.45\linewidth]{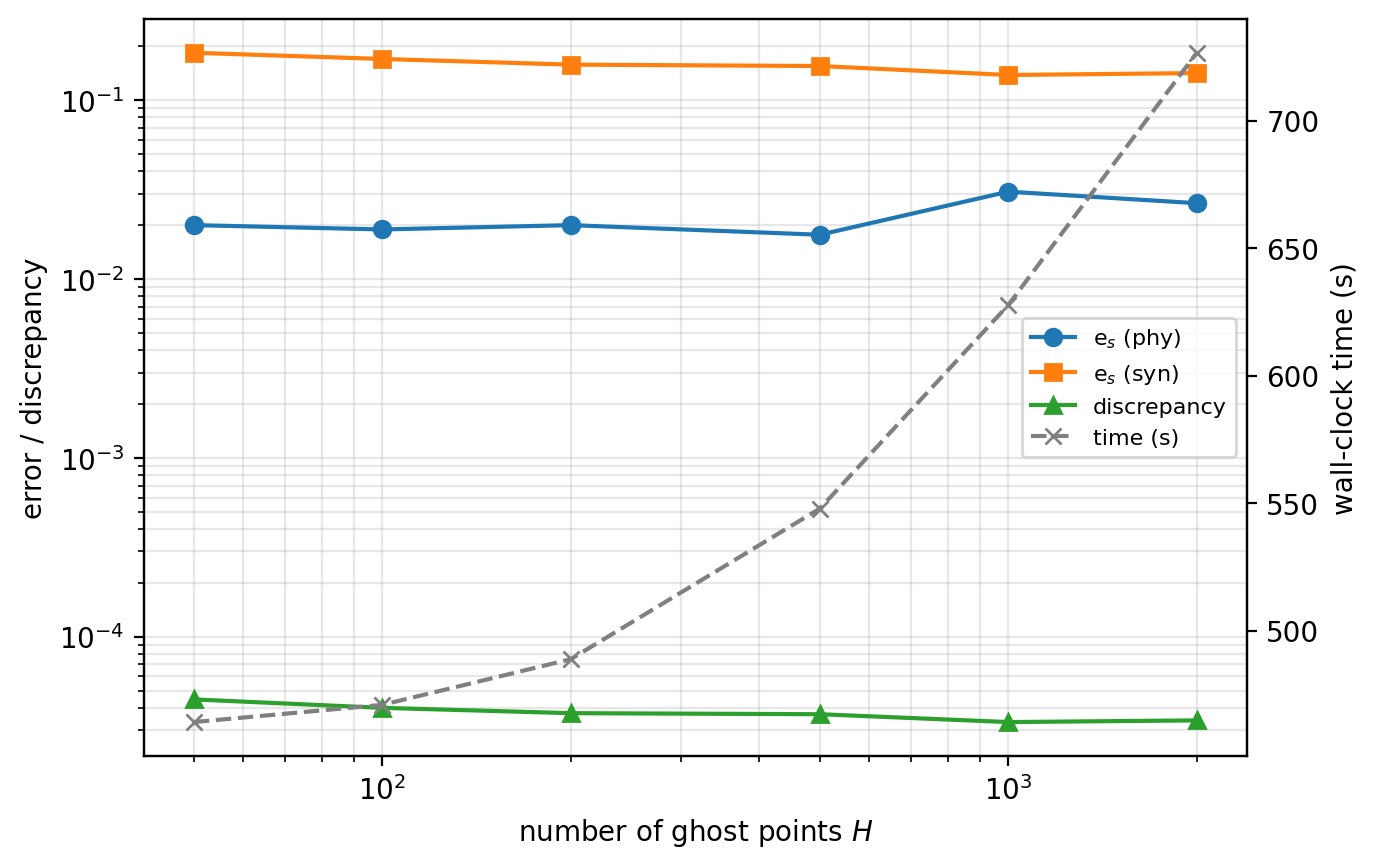}
    \includegraphics[width=0.45\linewidth]{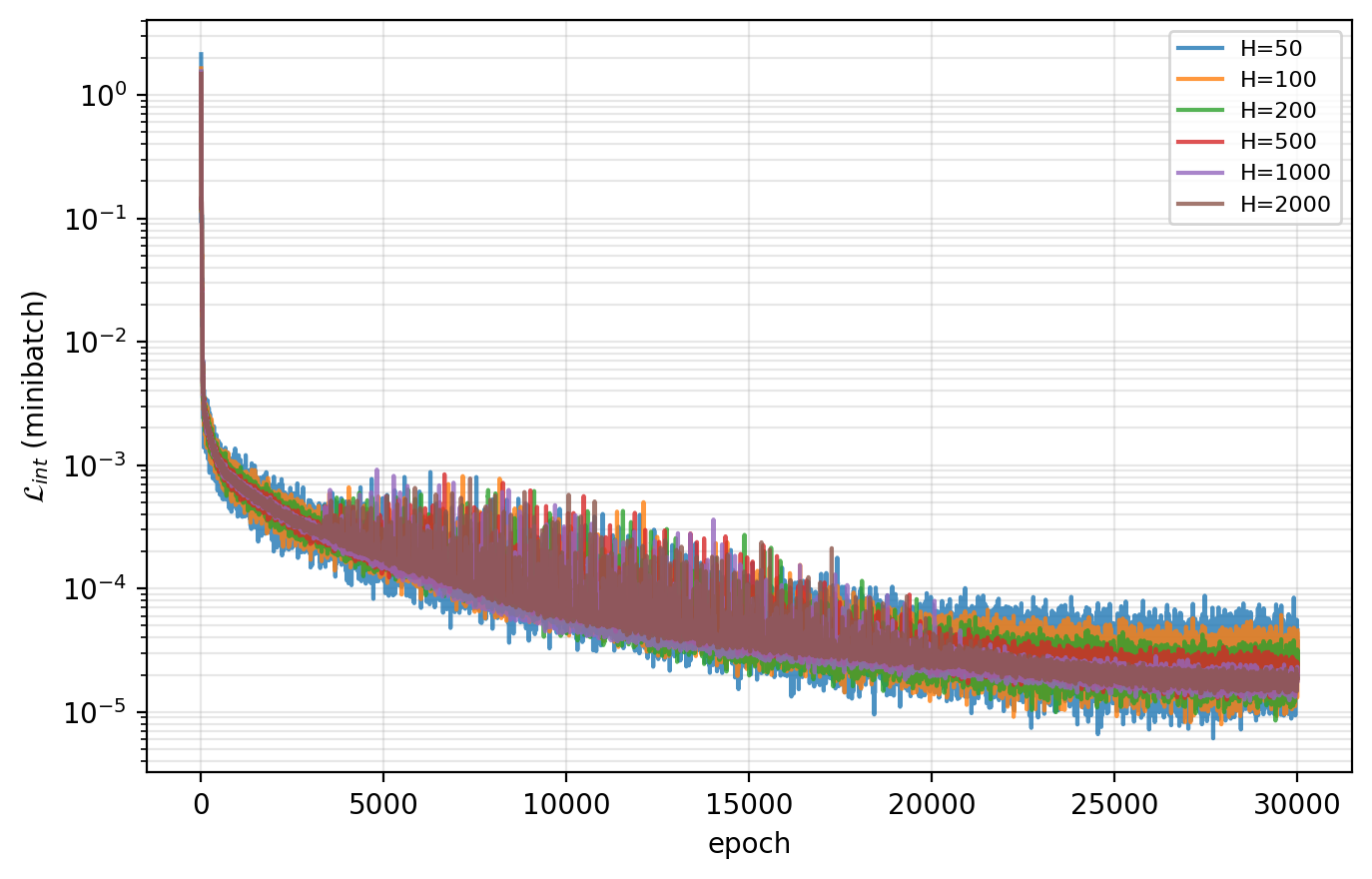}
    \caption{\textbf{Left:} Helmholtz sensitivity to the ghost-point count $H$. Physical and
    synthetic solution errors and the model--model discrepancy (left axis) against
    wall-clock time (right axis). Accuracy is flat in $H$ while cost grows. \textbf{Right:} Minibatch interaction loss $\mathcal{L}_{int}$ during training for
    each ghost-point count $H$. All values of $H$ converge to the same floor.}
    \label{fig:helm_ghost}
\end{figure}

\smallskip\noindent
\textbf{Ghost points.} Here, instead, we sweep the number of ghost points used in the interaction term
\[
H \in \{50, 100, 200, 500, 1000, 2000\}.\]
The results are summarized in Figure~\ref{fig:helm_ghost}.
On the left panel, we see that the physical solution error, the parameter error, and the model--model discrepancy
are essentially independent of $H$ ($\mathsf{e}^m_s(\mathrm{phy})\in[0.018,0.031]$,
$\mathsf{e}^m_p\in[0.012,0.028]$), while wall-clock time grows steadily from
$464$\,s to $727$\,s. Even $H=50$ ghost points, resampled per epoch, suffice to
enforce consistency between the two models; increasing $H$ buys no accuracy and
only adds cost. This validates the Random-Batch-style ghost sampling of
Section~\ref{sec:absstra} and justifies the small default $H=200$, although of course the ``correct'' number of ghost points depends on the complexity of the ground-truth model.

While the final errors are mostly independent on $H$, we point out that increasing the number of ghost points has an effect on the training dynamics. Indeed, in the right panel of Figure~\ref{fig:helm_ghost} we can observe that the interaction-loss histories oscillate less as $H$ increases, implying that a stronger interaction term tends to stabilize the HYCO training procedure.

\begin{figure}[h!]
    \centering
    \includegraphics[width=0.85\linewidth]{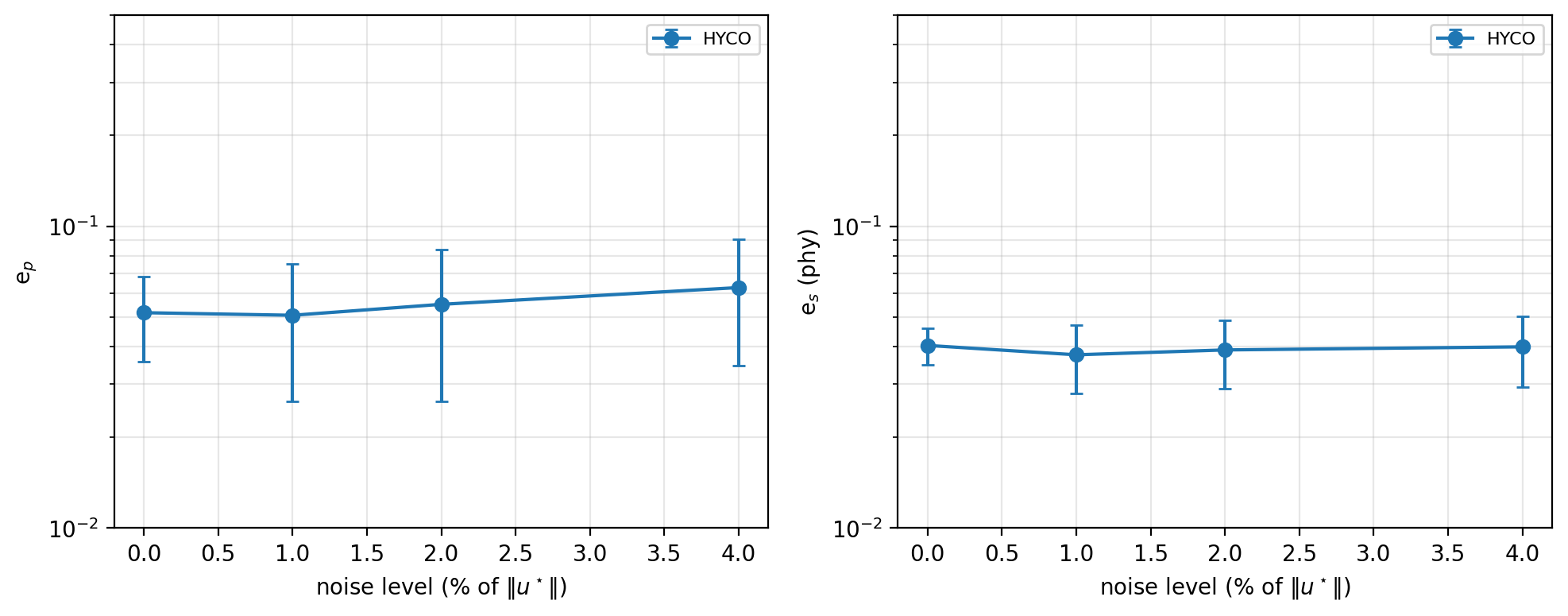}
    \caption{Helmholtz sensitivity to observational noise. Parameter error
    $\mathsf{e}^m_p$ (left) and physical solution error $\mathsf{e}^m_s$ (right) as
    the noise standard deviation grows to $4\%$ of $\|u^\star\|_{\mathrm{RMS}}$.}
    \label{fig:helm_noise}
\end{figure}

\smallskip\noindent
\textbf{Observational noise.} Finally, we add zero-mean Gaussian noise to the
sensor values, with standard deviation a prescribed fraction of
$\|u^\star\|_{\mathrm{RMS}}$, the root mean square of the reference solution over
the nodes of the fine mesh. Figure~\ref{fig:helm_noise} reports the recovered
errors at $0\%$, $1\%$, $2\%$, and $4\%$, across five different simulations with random seeds. Both $\mathsf{e}^m_p$ and
$\mathsf{e}^m_s$ stay within a narrow band across the whole range (parameter error
between $0.011$ and $0.029$), remaining at the same order of magnitude as the
noiseless case. Overall, parameter recovery is stable up to the $4\%$ noise level considered, showing a good robustness of the HYCO procedure to noise.

% =====================================================================
\section{A Dynamic Example: The LWR Traffic Model}
\label{sec:claw}

\noindent
We now turn to a first-order hyperbolic conservation law, a setting fundamentally different from the smooth problems considered so far: both the Gray-Scott pattern-formation system and the elliptic Helmholtz inversion admit smooth solutions, whereas here the solution develops genuine discontinuities. We consider the Lighthill-Whitham-Richards (LWR) model \cite{LWR}, the classical macroscopic description of vehicle density on a road, in which traffic evolves according to a scalar conservation law with a concave flux. Its solutions generically develop shocks, corresponding to a traffic jam forming, and rarefactions, corresponding to traffic dispersing, even from smooth initial data. This makes LWR a natural stress test for any PDE-learning method built around a smooth surrogate, and it is precisely why variants of it are used as a benchmark for physics-informed and operator-learning approaches to traffic-flow estimation and parameter identification from sensor data. We use it here as an inverse problem: identifying the model's physical parameters from partial density observations.

\subsection{Problem setup and HYCO components}

The density $\rho^\star(t,x)$ obeys the scalar conservation law
\begin{equation}
\label{eq:lwr}
\partial_t \rho + \partial_x\!\big(F(\rho;\Lambda^\star)\big) = 0,
\qquad (t,x) \in (0,T]\times(0,L),
\end{equation}
with the Greenshields flux
\[
F(\rho;\Lambda) = v_{\max}\,\rho\Big(1 - \frac{\rho}{\rho_{\max}}\Big),
\]
and a piecewise-constant initial condition mimicking a traffic jam,
\[
\rho_0(x) =
\begin{dcases}
0.9, & x \in [0.5, 1.5], \\
0.1, & \text{otherwise.}
\end{dcases}
\]
The unknown physical parameters are the maximal velocity and the maximal density,
\[
\Lambda^\star := (v_{\max}^\star, \rho_{\max}^\star) = (1,1).
\]
The reference solution $\rho^\star$ is computed on a fine finite-volume (FV) mesh of $2000$ cells, using a Godunov-type scheme with a Lax-Friedrichs numerical flux \cite{Godunov, LeVeque}, which resolves the shocks correctly. Observations are collected at $N_d = 100$ scattered space-time locations, restricted to the first quarter of the time interval ($t \in [0, 0.5]$ out of $T=2$) and the left half of the spatial domain ($x \in [0, 1.5]$ out of $[0,3]$), simulating sensors that cover only a portion of the road for a limited time. Recovering $\Lambda^\star$ therefore also requires extrapolating in time and space beyond the observation window.

\smallskip
\textbf{Goal:} {\it Given the partial density observations $\D$ and the model form \eqref{eq:lwr}, identify the parameters $(v_{\max}, \rho_{\max})$ that best explain the observed traffic evolution.}

\smallskip
\noindent
The HYCO components are:
\begin{itemize}
\item[$\diamond$] The \textbf{physical model} is the FV solver of \eqref{eq:lwr}, parametrized by $\Lambda = (v_{\max}, \rho_{\max})$, run on a coarse mesh of $50$ cells, forty times coarser than the $2000$-cell mesh used to generate the data. This coarse solver still captures the shock structure while keeping the forward solves cheap. Its data-fitting loss is \eqref{loss_phy}, specialized to the $N_d$ observation points.
\item[$\diamond$] The \textbf{synthetic model} is a feed-forward neural network $\rho_{syn}:(0,T)\times(0,L) \to \R_{\geq 0}$ mapping a space-time point to the density, with $2$ hidden layers of width $128$, ReLU activation, and residual connections, as in \eqref{eq:nn} (narrower hidden width than in Section~\ref{sec:helm}), and output constrained to be nonnegative. Its loss mirrors \eqref{loss_syn}, evaluated at the same observation points.
\item[$\diamond$] The \textbf{interaction term} penalizes the discrepancy of the two models on $H$ ghost points $\{(t_h^G, x_h^G)\}$ sampled in the space-time domain,
\[
\L_{int}(\Theta,\Lambda) = \frac{1}{H}\sum_{h=1}^H \|\rho_{syn}(t_h^G, x_h^G) - \rho_{phy}(t_h^G, x_h^G)\|^2 .
\]
\end{itemize}
As before, the physical model minimizes $\L_1 = \beta\L_{phy} + \L_{int}$ over $\Lambda$ and the synthetic model minimizes $\L_2 = \alpha\L_{syn} + \L_{int}$ over $\Theta$, alternating updates. We set $\alpha = 1$ and $\beta = 0$: the physical model receives no direct data term and is driven by the interaction loss alone.

\smallskip
\noindent
\textbf{Numerical parameters.} All models are trained with the Adam optimizer, with a maximum budget of $20{,}000$ epochs. In HYCO, the physical and synthetic models use cosine-decayed learning rates starting at $1\times10^{-3}$ and $5\times10^{-4}$, respectively, decaying to $1\%$ of their initial value, with $H = 2000$ ghost points per epoch. The PINN baseline uses the same network width, $\tanh$ activation, a cosine-decayed learning rate starting at $1\times10^{-4}$, a residual weight $\gamma = 5$, and $2000$ collocation points. Training stops once the stabilization criterion \eqref{stopping} is met, with window $Z=1000$ and tolerance $\varepsilon=10^{-4}$ for HYCO and PINN, and $Z=5$, $\varepsilon=10^{-4}$ for EKI. All results reported below are averaged over five random seeds $\{42,33,17,23,2\}$ and given as mean $\pm$ standard deviation.

\subsection{Basic experiment}

\begin{figure}[h!]
    \centering
    % PLACEHOLDER
    \includegraphics[width=\linewidth]{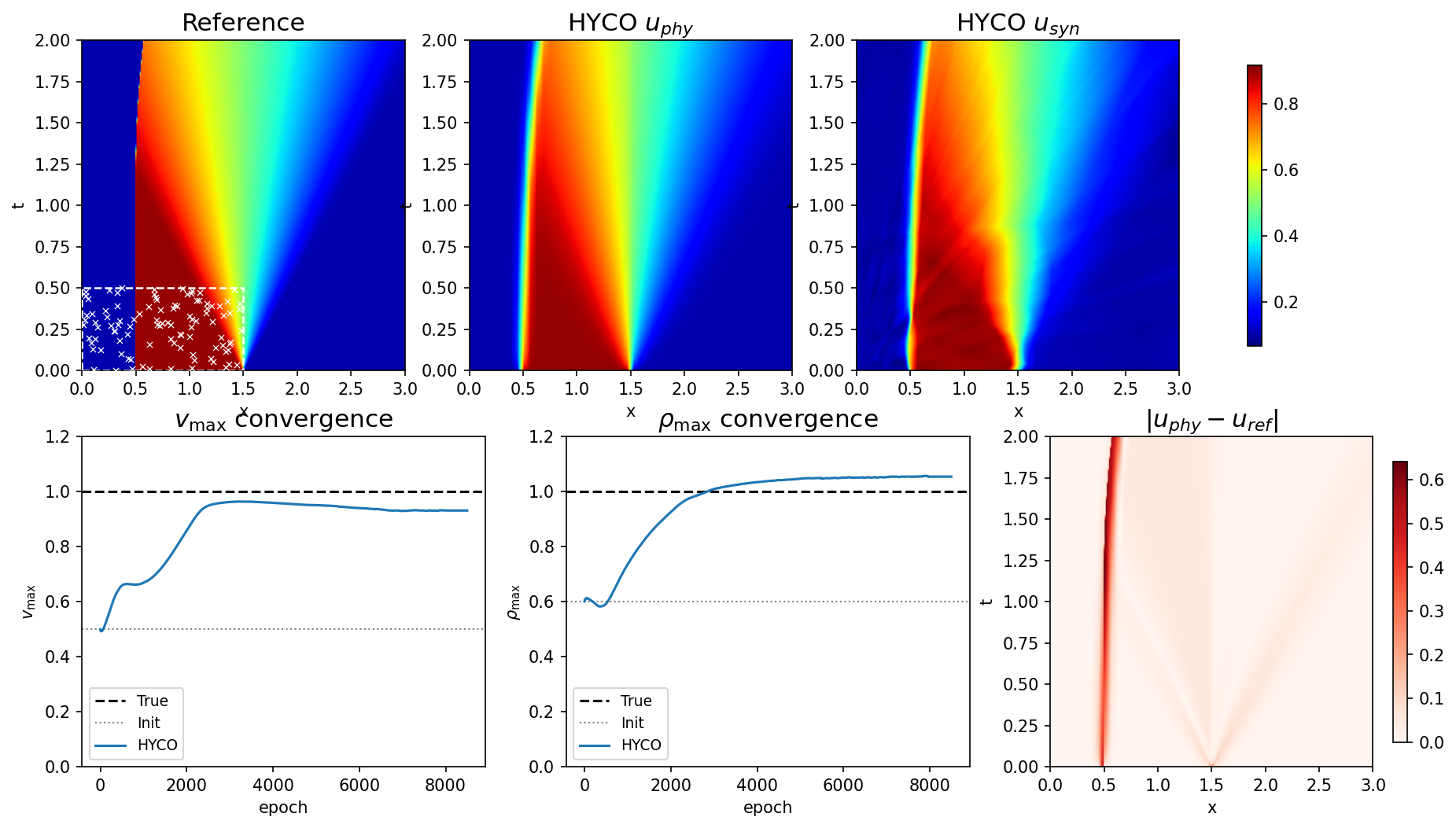}
    \caption{LWR basic experiment. \textbf{Row 1:} space-time density fields for the reference $\rho^\star$, HYCO's physical prediction $\rho_{phy}$, and HYCO's synthetic prediction $\rho_{syn}$; the dashed box and white $\times$ markers in the Reference panel mark the observation window and sensor locations. \textbf{Row 2:} convergence of $v_{\max}$ and $\rho_{\max}$ over training (dashed: ground truth, dotted: initial guess), and the pointwise error $|\rho_{phy}-\rho^\star|$ in space-time. Training stops once the stabilization criterion \eqref{stopping} is satisfied, which is why the convergence curves end before the nominal epoch budget. Fields shown for seed $42$.}
    \label{fig:claw_basic}
\end{figure}

Figure~\ref{fig:claw_basic} shows a single HYCO run: the space-time density fields for the reference, the physical predictor, and the synthetic surrogate, together with density profiles at two representative times, one inside the observation window and one in the extrapolation region, and the convergence of $v_{\max}$ and $\rho_{\max}$ over training. HYCO recovers both parameters and reproduces the shock propagation well beyond the observed space-time region.

\subsection{Ablation and comparison}

\noindent
As for Helmholtz, we now assess how much of HYCO's performance is due to the cooperative coupling, and how it compares to existing inversion strategies, here in a setting where the presence of shocks makes both questions considerably harder.

\smallskip
\noindent
\textbf{Ablation.} We disable the interaction term ($\L_{int}=0$). Since $\beta = 0$, this removes the only signal driving the physical model: $\Lambda$ is never updated and remains at its initial value for the whole training run. The ablation is therefore only meaningful for the synthetic model, which keeps its own data term $\alpha\L_{syn}$ but loses the physics-informed consistency the interaction loss otherwise provides outside the observed region. Table~\ref{tab:claw_comp} (row HYCO-abl) reports the resulting solution error $\mathsf{e}^m_s$; its large value reflects that, without the interaction with the physical model, the only incentive to the synthetic model is to fit the data in the small collection region.

\begin{figure}[h!]
    \centering
    % PLACEHOLDER
    \includegraphics[width=\linewidth]{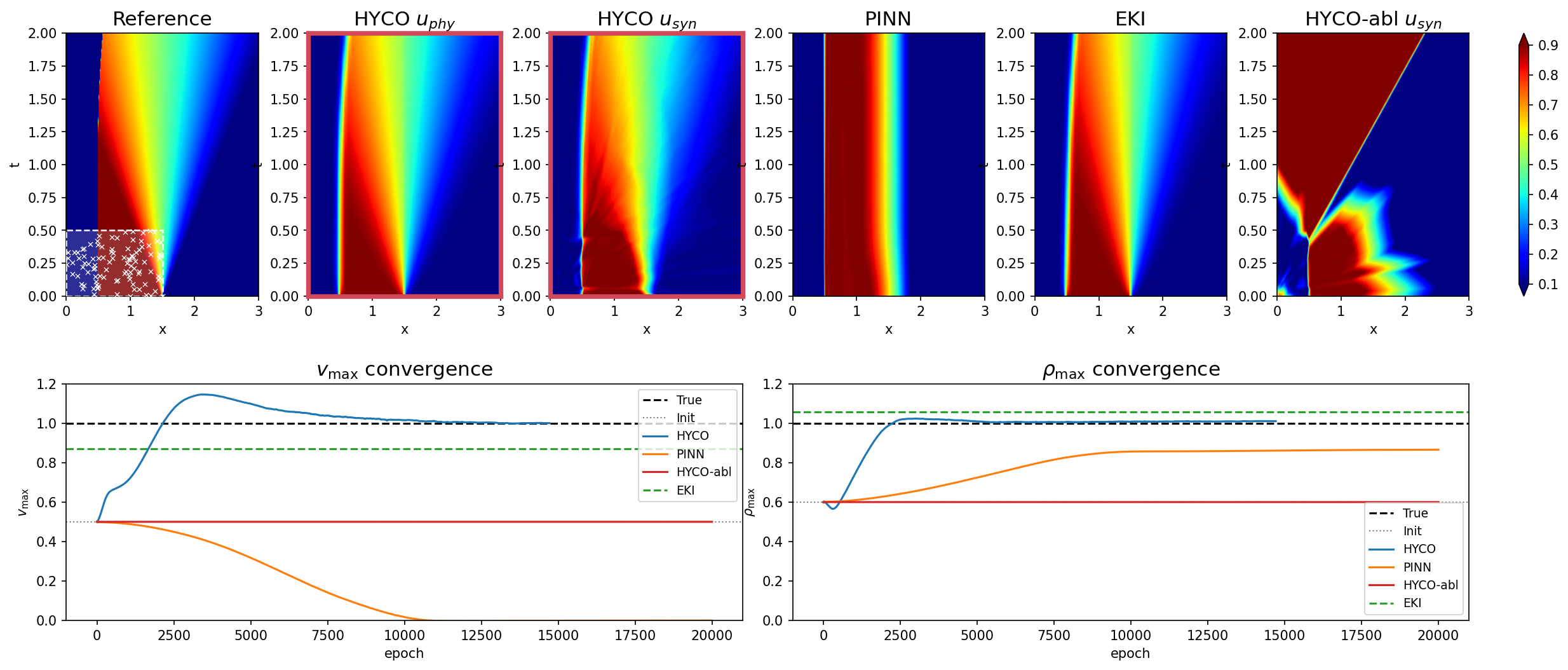}
    \caption{LWR comparison. \textbf{Row 1:} space-time density fields for the reference $\rho^\star$ (observation window and sensor locations marked as in Figure~\ref{fig:claw_basic}), HYCO's physical and synthetic predictions, PINN, EKI, and the ablated HYCO-abl synthetic model. \textbf{Row 2:} convergence of $v_{\max}$ and $\rho_{\max}$ over training for HYCO and PINN (dashed: ground truth, dotted: initial guess); EKI is shown as a horizontal dashed line at its converged ensemble estimate, since it has no comparable notion of training epoch. The HYCO and PINN curves end at different epochs because each stops once its own stabilization criterion \eqref{stopping} is satisfied, rather than at a fixed epoch budget. Fields shown for seed $42$.}
    \label{fig:claw_comp}
\end{figure}

\smallskip
\noindent
\textbf{Competitors.} A shock-forming conservation law is a harder setting for both baseline families used in Section~\ref{sec:helm}, and this shapes which variants we include here. As in Section~\ref{sec:helm}, we exclude Neural Operator baselines for the same reason: this is a single-instance inverse problem, not the multi-instance operator-learning setting they are designed for.
\begin{enumerate}[label=(\roman*)]
\item \textbf{PINN} \cite{RPK}: a $\tanh$ network trained with a data loss plus the conservation-law residual, with $\Lambda$ learned jointly. We stress that this is not a natural method for conservation laws: a smooth $\tanh$ ansatz cannot represent discontinuous shocks by construction, so the PINN is included as a reference point rather than a competitive baseline, and to make this failure mode explicit. We omit the self-adaptive variant (SA-PINN) used for Helmholtz: reweighting collocation points addresses where the residual is large, not the underlying inability of a smooth ansatz to represent a discontinuity, so it would not fix PINN's failure mode here.
\item \textbf{Ensemble Kalman Inversion (EKI)} \cite{EKI}: a gradient-free ensemble method tailored to solver-based inversion. An ensemble of $J = 50$ parameter vectors is initialized around $\Lambda_0$ and, at each iteration, propagated through the FV solver and updated via the perturbed-observation Kalman gain; the estimate is the ensemble mean. EKI maintains no solution surrogate, so its solution error $\mathsf{e}^m_s$ is evaluated by running the forward solver with the recovered $\hat\Lambda$. We use EKI here in place of the Gauss-Newton baseline used for Helmholtz because it is gradient-free: Gauss-Newton requires differentiating through the solver, which is unreliable close to a shock, whereas EKI only needs forward solver evaluations.
\end{enumerate}
Table~\ref{tab:claw_comp} collects the results, with recovered densities and convergence curves in Figure~\ref{fig:claw_comp}.

\begin{table}[h!]
\centering
\small
\renewcommand{\arraystretch}{1.2}
\setlength{\tabcolsep}{4pt}
\begin{tabular}{|c|c|c|c|c|c|c|}
\hline
\textbf{Method} & $\hat v_{\max}$ & $\hat \rho_{\max}$ & $\mathsf{e}^m_p$ & $\mathsf{e}^m_s$ & $\mathsf{e}^m_d$ & \textbf{ms/epoch\,$\vert$\,iter} \\
\Xhline{1.2pt}
PINN & \makecell{0.0020 \\ {\scriptsize $\pm0.0047$}} & \makecell{0.8415 \\ {\scriptsize $\pm0.0521$}} & \makecell{0.7155 \\ {\scriptsize $\pm0.0032$}} & \makecell{1.1586 \\ {\scriptsize $\pm0.3360$}} & \makecell{1.2560e-03 \\ {\scriptsize $\pm5.65$e-04}} & 7.4 \\
\hline
EKI & \makecell{0.8437 \\ {\scriptsize $\pm0.0941$}} & \makecell{1.0569 \\ {\scriptsize $\pm0.0417$}} & \makecell{0.1304 \\ {\scriptsize $\pm0.0462$}} & \makecell{0.1597 \\ {\scriptsize $\pm0.0476$}} & \makecell{5.1393e-03 \\ {\scriptsize $\pm3.49$e-03}} & 112.8$^\dagger$ \\
\hline
HYCO-abl ($\L_{int}=0$) & -- & -- & -- & \makecell{2.0681 \\ {\scriptsize $\pm0.4520$}} & \makecell{7.6729e-05 \\ {\scriptsize $\pm1.36$e-04}} & 11.3 \\
\hline
HYCO & \makecell{1.0025 \\ {\scriptsize $\pm0.1069$}} & \makecell{1.0229 \\ {\scriptsize $\pm0.0237$}} & \makecell{\textbf{0.0595} \\ {\scriptsize $\pm0.0522$}} & \makecell{\textbf{0.1378} \\ {\scriptsize $\pm0.0549$}} & \makecell{\textbf{6.9020e-05} \\ {\scriptsize $\pm4.39$e-05}} & 13.0 \\
\hline
\end{tabular}
\caption{LWR: ablation and comparison. Recovered parameters $\hat v_{\max}, \hat \rho_{\max}$, parameter error $\mathsf{e}^m_p$ \eqref{eq:normeu}, solution error $\mathsf{e}^m_s$ \eqref{eq:normL2}, and data discrepancy $\mathsf{e}^m_d$ \eqref{eq:dataerr}. Ground truth $(v_{\max}^\star,\rho_{\max}^\star)=(1,1)$. Best accuracy in bold. Error entries are mean $\pm$ standard deviation over five random seeds $\{42,33,17,23,2\}$. The last column is the per-step wall-clock cost on common hardware (single $1000$-epoch run, seed $42$): milliseconds per epoch for the gradient-trained models, and per iteration for EKI. $^\dagger$EKI is not epoch-based and converges in only $30$ iterations ($3.4$\,s total), so its per-iteration cost is not directly comparable to an epoch.}
\label{tab:claw_comp}
\end{table}

\noindent
\textbf{Conclusions.}
\begin{itemize}
\item[$\diamond$] The PINN fails to represent the shock and, consequently, misidentifies the parameters, as anticipated by the smoothness of its ansatz. We report it nonetheless to make explicit why a shock-aware comparator is needed.
\item[$\diamond$] EKI, operating directly on the shock-resolving FV solver, is a strong classical competitor; nonetheless, averaged over the five seeds, HYCO recovers the parameters more accurately ($\mathsf{e}^m_p = 0.0595$, against EKI's $0.1304$) and also attains the smallest solution error ($\mathsf{e}^m_s = 0.1378$, against EKI's $0.1597$). In addition, HYCO provides a synthetic surrogate valid across the whole space-time domain, including the unobserved region, whereas EKI returns only the recovered parameters. The comparison is not one-sided, however: EKI is by far the cheaper method in absolute terms, converging in $30$ iterations and about $3.4$\,s, well below the cost of a HYCO run at $13.0$\,ms per epoch over several thousand epochs. HYCO's advantage here is recovery accuracy and the additional surrogate, not computational cost.
\item[$\diamond$] Here the per-epoch picture is reversed: HYCO is about $1.75$ times more expensive per epoch than the PINN ($13.0$ vs $7.4$\,ms), since every physical update differentiates through a full coarse finite-volume solve whereas the PINN only evaluates a first-order residual. The cooperative coupling itself adds only about $15\%$ over the decoupled baseline ($11.3$\,ms), and the PINN's cheaper step is moot, as its smooth ansatz cannot represent the shock. EKI is likewise very cheap in absolute terms, converging in $30$ iterations and about $3.4$\,s; but, as for Gauss-Newton, its accuracy plateaus well before the iteration cap (Figure~\ref{fig:claw_comp}), so the residual gap to HYCO on both $\mathsf{e}^m_p$ and $\mathsf{e}^m_s$ is not one that additional iterations would remove. HYCO's advantage on this problem is thus convergence quality rather than a cheaper step.
\end{itemize}

\section{A Toy Theoretical Setting}
\label{sec:games}

\noindent
The non-convexity of the HYCO loss functionals, inherited from the nonlinear dependence of the synthetic model on its parameters, puts a full convergence analysis out of reach. Here we settle for a more modest but rigorous statement: in a simplified, convex regime, the HYCO optimization game admits a Nash equilibrium. Concretely, we pair a linear (Poisson) physical model with a mean-field relaxation of a shallow network, and read the alternating scheme as a two-player game in which the physical and synthetic models each minimize their own objective while coupled through the shared consistency term. We stress that here we only establish the existence of an equilibrium; we do not address whether the alternating training algorithm used in practice converges to this equilibrium; this is left for future work.

Formally, we introduce two players. The first is the physical model, which minimizes
\[
\L_1(\Theta, \Lambda) := \beta \L_{phy}(\Lambda) + \L_{int}(\Theta, \Lambda)
\]
with respect to its own parameters $\Lambda$, keeping $\Theta$ fixed, subject to $\u_{phy}$ solving the physical model. The second is the synthetic model, which minimizes
\[
\L_2(\Theta, \Lambda) := \alpha \L_{syn}(\Theta) + \L_{int}(\Theta, \Lambda)
\]
with respect to $\Theta$, keeping $\Lambda$ fixed.

\begin{definition}
\label{def:eq}
A pair $(\Theta^\star, \Lambda^\star)$ is a \emph{Nash equilibrium} if
\[
\L_1(\Theta^\star, \Lambda^\star) \leq \L_1(\Theta^\star, \Lambda) \quad \text{and} \quad 
\L_2(\Theta^\star, \Lambda^\star) \leq \L_2(\Theta, \Lambda^\star)
\]
for all admissible $\Theta$ and $\Lambda$.
\end{definition}

\smallskip
\noindent
\textbf{A potential-game viewpoint.} The two objectives are not arbitrary. Both $\L_1$ and $\L_2$ derive from the single joint functional
\[
\L(\Theta,\Lambda) := \alpha\L_{syn}(\Theta) + \beta\L_{phy}(\Lambda) + \L_{int}(\Theta,\Lambda)
\]
by freezing one block of variables: since $\alpha\L_{syn}$ is independent of $\Lambda$ and $\beta\L_{phy}$ of $\Theta$, one has $\partial_\Lambda \L = \partial_\Lambda \L_1$ and $\partial_\Theta \L = \partial_\Theta \L_2$. The game is thus an \emph{exact potential game} with potential $\L$.

Two distinct objects should be kept apart at this point. In the \emph{idealized} scheme, the full interaction integral \eqref{loss_int} is used and each block is minimized exactly; this is genuine block-coordinate descent on $\L$, so each block minimization decreases $\L$ monotonically, and any fixed point of the alternating iteration is a Nash equilibrium in the sense of Definition~\ref{def:eq}. The \emph{implemented} algorithm of Sections~\ref{sec:grayscott},~\ref{sec:helm} and~\ref{sec:claw} is an inexact stochastic approximation of it: each block receives a single gradient step rather than an exact minimization, and the interaction integral is replaced by an empirical average over ghost points resampled at every epoch, so that the sampled potential itself changes from one epoch to the next. Accordingly, the potential-game structure is offered as a structural motivation for the alternating procedure, and does not by itself imply monotonicity of $\L$ or convergence of the algorithm used in the experiments.

Before diving into the simplified HYCO formulation, we report an existence result for two-player convex games. The proposition is standard, being a direct combination of classical facts \cite{Fan1952, Debreu1952, Glick}, but we give a short proof since we could not locate it stated under exactly the hypotheses we need.

\begin{proposition}
\label{prop:nash}
Let $S_1, S_2$ be nonempty, convex, and compact metrizable subsets of two real locally convex Hausdorff topological vector spaces, and let $J_1, J_2 : S_1 \times S_2 \to \R$ be two payoffs such that, for each $i \in \{1,2\}$:
\begin{enumerate}[label=(\roman*)]
\item the map $s_i \mapsto J_i(s_1, s_2)$ is convex;
\item $J_i$ is sequentially lower semicontinuous on $S_1 \times S_2$ and continuous in the variable of the other player.
\end{enumerate}
Then the game admits a Nash equilibrium, i.e. a pair $(s_1^\star, s_2^\star) \in S_1 \times S_2$ with
\[
J_1(s_1^\star, s_2^\star) \le J_1(s_1, s_2^\star) \quad\text{and}\quad J_2(s_1^\star, s_2^\star) \le J_2(s_1^\star, s_2)
\]
for all $s_1 \in S_1$, $s_2 \in S_2$.
\end{proposition}

\begin{proof}
For each fixed $s_2$, the best-response set $\mathrm{BR}_1(s_2) := \argmin_{s_1 \in S_1} J_1(s_1, s_2)$ is nonempty, since a sequentially lower semicontinuous function attains its minimum on the compact metrizable set $S_1$, and convex by (i); the same holds for $\mathrm{BR}_2$. If now $s_1^n \in \mathrm{BR}_1(s_2^n)$ and $(s_1^n, s_2^n) \to (s_1, s_2)$, letting $n \to \infty$ in $J_1(s_1^n, s_2^n) \le J_1(t, s_2^n)$ yields $J_1(s_1, s_2) \le J_1(t, s_2)$ for every $t \in S_1$, by lower semicontinuity on the left and continuity in the opponent's variable on the right; hence $s_1 \in \mathrm{BR}_1(s_2)$. By metrizability, the correspondence $(s_1, s_2) \mapsto \mathrm{BR}_1(s_2) \times \mathrm{BR}_2(s_1)$ therefore has closed graph, with nonempty convex values on the convex compact set $S_1 \times S_2$. The Kakutani--Fan--Glicksberg fixed-point theorem \cite{Glick} yields a fixed point, which is a Nash equilibrium.
\end{proof}

We now instantiate Proposition~\ref{prop:nash} in the HYCO setting. The physical ansatz is the Poisson equation with homogeneous Dirichlet boundary conditions,
\begin{equation}
\begin{cases}
    -\Delta \u_{phy}(\x) = f(\x), & \x \in \Omega, \\
    \u_{phy}(\x) = 0, & \x \in \partial \Omega,
\end{cases}
\end{equation}
where $\Omega \subset \R^d$, $d \le 3$, is a bounded domain with $C^{1,1}$ boundary. The source term $f \in L^2(\Omega)$ is the parameter to be optimized, constrained to the ball $\{\|f\|_{L^2} \le C\}$ for a fixed $C > 0$. The solution operator $(-\Delta)^{-1}$ is compact and continuous from $L^2(\Omega)$ to $L^2(\Omega)$. Given observations $\D = \{\u^\star(\x_i)\}_{i=1}^M$, the physical loss is
\[
\L_{phy}(f) := \frac{1}{M} \sum_{i=1}^M \| ((-\Delta)^{-1}f)(\x_i) - \u^\star(\x_i) \|^2.
\]
The pointwise evaluations are well defined: elliptic regularity gives $\u_{phy} \in H^2(\Omega)$, and since $d \le 3$ one has $H^2(\Omega) \hookrightarrow C^0(\overline\Omega)$.

For the synthetic model we use the convex relaxation of a shallow network with $P$ neurons,
\begin{equation}
\label{nn_delta}
\u_{syn}(\x) = \sum_{j=1}^P w_j \sigma(\langle a_j, \x\rangle + b_j),
\end{equation}
with ReLU activation $\sigma$. The associated sparsity-regularized empirical loss
\begin{equation}
\label{sparse_reg}
\min_{\Theta} \lambda \sum_{j=1}^P |w_j| + \frac{1}{M} \sum_{i=1}^M \left\| \sum_{j=1}^P w_j \sigma(\langle a_j, \x_i\rangle + b_j) - \u^\star(\x_i) \right\|^2
\end{equation}
is nonconvex. Following \cite{KANG}, we relax the parameter space to finite measures: letting $K \subset \R^{d+1}$ be a compact set containing all $(a,b)$ parameters, we replace the finite sum by
\begin{equation}
\label{nn_relax}
\Tilde{\u}_\mu(\x) := \int_{K} \sigma(\langle a, \x\rangle + b) \, d\mu(a, b),
\qquad \mu \in \mathcal{M}(K),
\end{equation}
so that the relaxed synthetic loss
\begin{equation}
\label{reg_loss}
\L_{syn}(\mu) := \lambda \|\mu\|_{TV} + \frac{1}{M} \sum_{i=1}^M \left\| \Tilde{\u}_\mu(\x_i) - \u^\star(\x_i) \right\|^2
\end{equation}
is convex in $\mu$. The coupling is the $L^2$ distance between the two predictions,
\[
\L_{int}(f, \mu) := \|(-\Delta)^{-1}f - \Tilde{\u}_\mu\|^2_{L^2(\Omega)},
\]
which is well defined since $\mu$ is finite and $\sigma$ is locally bounded. The two players then minimize
\[
\L_1(f, \mu) := \frac{\beta}{M} \sum_{i=1}^M \| ((-\Delta)^{-1}f)(\x_i) - \u^\star(\x_i) \|^2 + \|(-\Delta)^{-1}f - \Tilde{\u}_\mu\|^2_{L^2(\Omega)},
\]
\[
\L_2(f, \mu) := \lambda \|\mu\|_{TV} + \frac{\alpha}{M} \sum_{i=1}^M \left\| \Tilde{\u}_\mu(\x_i) - \u^\star(\x_i) \right\|^2 + \|(-\Delta)^{-1}f - \Tilde{\u}_\mu\|^2_{L^2(\Omega)},
\]
over the admissible sets
\[
\mathcal{B}_C := \{f \in L^2(\Omega): \|f\|_{L^2} \le C\}, \qquad
\mathcal{A}_R := \{\mu \in \mathcal{M}(K): \|\mu\|_{TV} \le R\},
\]
for a fixed $R > 0$.
We then have the following theorem.

\begin{theorem}
\label{thm:hyco-nash}
The relaxed HYCO game with payoffs $\L_1, \L_2$ and strategy sets $\mathcal{B}_C, \mathcal{A}_R$ admits at least one Nash equilibrium in the sense of Definition~\ref{def:eq}.
\end{theorem}

\begin{proof}
The claim follows from Proposition~\ref{prop:nash} once its hypotheses are checked. We equip $\mathcal{B}_C$ with the weak topology of $L^2(\Omega)$ and $\mathcal{A}_R$ with the weak-$*$ topology of $\mathcal{M}(K) = C(K)^*$; both are metrizable on these bounded sets, since $L^2(\Omega)$ is separable and so is $C(K)$, $K$ being compact.

First, we observe that both strategy sets are then admissible. Indeed, the set $\mathcal{B}_C$ is a closed bounded convex subset of a Hilbert space, hence weakly compact and convex; $\mathcal{A}_R$ is a closed ball in a dual space, hence convex and, by the Banach--Alaoglu theorem, weak-$*$ compact. Convexity of the payoffs is equally straightforward: each is a sum of the total-variation norm, which is convex, and squared norms of quantities depending affinely on the linear maps $f \mapsto (-\Delta)^{-1}f$ and $\mu \mapsto \Tilde{\u}_\mu$, so $\L_1$ is convex in $f$ and $\L_2$ is convex in $\mu$.

The only delicate point is continuity of the coupling in each variable. On the physical side, elliptic regularity makes $(-\Delta)^{-1}$ bounded from $L^2(\Omega)$ into $H^2(\Omega)$, while for $d \le 3$ the embedding $H^2(\Omega) \Subset C^0(\overline\Omega)$ \cite[\S5.6.3, \S6.3]{Evans} is compact; the composition $L^2(\Omega) \to C^0(\overline\Omega)$ is therefore compact and turns $f_n \rightharpoonup f$ into uniform convergence of the states. On the synthetic side, $(a,b) \mapsto \sigma(\langle a, \x \rangle + b)$ is continuous on the compact $K$, so weak-$*$ convergence gives $\Tilde{\u}_{\mu_n} \to \Tilde{\u}_\mu$ pointwise, uniformly bounded thanks to $\|\mu_n\|_{TV} \le R$; dominated convergence then yields strong $L^2$ convergence; hence the terms depending on $\mu$ are continuous, while $\|\mu\|_{TV}$ is weak-$*$ lower semicontinuous. Each payoff is therefore continuous in the opponent's variable and jointly lower semicontinuous, and Proposition~\ref{prop:nash} gives the equilibrium.
\end{proof}

The numerical experiments in Sections~\ref{sec:grayscott}, \ref{sec:helm}, and \ref{sec:claw} show that HYCO performs well empirically far outside this convex, relaxed regime. Extending the guarantees of this section to the nonconvex setting used in practice, and establishing convergence of the alternating algorithm itself, remain open directions.

\section{Conclusions and perspectives}
\label{sec:conclusions}

\noindent
In this work, we introduced HYCO, a hybrid modeling strategy that integrates physical and synthetic models by optimizing them in an alternating fashion while encouraging alignment between their predictions. HYCO treats both models as active participants in the learning process, allowing mutual influence that serves as a powerful form of regularization.

Through numerical experiments on both steady and time-dependent PDEs, we showed that HYCO performs well in a wide range of tasks, from learning the solution of an unknown physical system, to parameter identification. These gains are particularly evident when observations are spatially localized, as the interplay between the models enables the reconstruction of missing information.

Looking ahead, we outline different research directions.
\begin{enumerate}
    \item We have validated the general empirical performances of HYCO on some classical examples. However, HYCO is a very general and flexible architecture. Future works will focus on applying the technique to several scenarios:
    \begin{itemize}
    \item[$\diamond$] In multiphysics scenarios (such as thermo-mechanical systems or electrochemical devices), each model may be associated with a different physical process and trained on distinct data. By properly designing the interaction term, the HYCO strategy can serve as a principled way to integrate these heterogeneous sources into a unified predictive framework. 
    \item[$\diamond$] Similarly, HYCO can be applied to problems involving multiple spatial domains. By tailoring the interaction term to act only at the interfaces, the strategy naturally enforces continuity between the models' predictions, making it suitable for domain decomposition methods or coupled simulations. 
    \item[$\diamond$] More broadly, a key step for future work will be to validate HYCO on more complex and large datasets arising from real-world applications (climate modeling, robotics, etc.), where physical fidelity and data scarcity coexist. In this context, it will be interesting to investigate the performances of HYCO using other popular machine learning architectures for the synthetic model. Transformers and Neural Operators are an obvious first choice, in light of the important results recently obtained by these architecture in data-driven learning of physical systems \cite{NO, TRA}.
    \end{itemize}
    \item Furthermore, while we touched upon the potential for privacy-aware versions of HYCO, many aspects remain to be explored. One natural direction is to investigate how the method performs when each model is trained on its own private dataset without sharing inputs. Privacy could also be enhanced through pre-processing techniques, such as Gaussian smoothing of the training data. The effectiveness and convergence of HYCO in these settings is not yet known and requires investigation. Moreover, as prior work in Federated Learning has shown that privacy breaches are possible even without data sharing \cite{SongWangZuazua}, it is critical to assess whether similar vulnerabilities might exist in our framework.
    \item A last challenge is to provide a mathematical understanding of HYCO's performance. Although neural networks complicate the theoretical analysis, it is important to note that the synthetic model is not restricted to being a machine learning architecture. This might allow, for instance, to obtain analytical results in simplified settings, providing some interpretability. In this direction, \cite{CoercivityGap} studies a related state-level stability phenomenon for neural PDE solvers: even when the neural parameters themselves fail to converge, the corresponding state functions may still converge strongly to the true solution, a principle shown to extend to HYCO-type hybrid methods. Incorporating this perspective into a convergence analysis for the alternating training algorithm itself is a natural next step.

    Besides, we also highlighted that a game-theoretic perspective offers a promising theoretical background. While this field is vast and complex, further exploration may uncover useful insights or training algorithms. For example, it would be interesting to investigate whether one might obtain results similar to those in \cite{KANG}, where it is shown that there is  no relaxation gap in the training of neural networks. In our context, it is natural to wonder whether the same is true for the Nash equilibrium described in Section~\ref{sec:games}.
\end{enumerate}

\section*{Acknowledgments}
\noindent
E.Z. has been funded by the project CoDeFeL of the European Research Council (ERC) under the European Union's Horizon Europe research and innovation program (grant agreement No. 101096251), and by the Alexander von Humboldt-Professorship program. This material is based upon work supported by the Air Force Office of Scientific Research under award number FA8655-22-1-7012, in the amount of \$160,907, grant PID2023-146872OB-I00 of the AEI (Spain), the COST Actions CA24136 -- InterCoML and CA24122 -- mSPACE, and the SURE-AI Center grant 357482, Research Council of Norway.

%%%% References %%%%
\bibliographystyle{abbrv}
\bibliography{bib}

\appendix
\appendix

\section{A Practical Guide to Applying HYCO}
\label{app:guide}

\noindent
This appendix collects the practical steps for applying HYCO to a new problem. We
separate the elements that define the framework from the implementation choices
made in our experiments, then give a recipe, guidance on the hyperparameters, and
remarks on when HYCO is most useful.

\subsection{Mandatory ingredients versus implementation choices}
The HYCO framework, as formulated in \eqref{main_minimization} and
Algorithm~\ref{alg:hyco}, requires only four elements: (i) a physical model that
maps parameters $\Lambda$ to a prediction $\u_{phy}$; (ii) a synthetic model that
maps parameters $\Theta$ to a prediction $\u_{syn}$; (iii) the three losses
$\L_{phy}, \L_{syn}, \L_{int}$, of which the interaction term $\L_{int}$ is the
only one coupling the two models; and (iv) alternating minimization of
$\L_1 = \beta\L_{phy} + \L_{int}$ and $\L_2 = \alpha\L_{syn} + \L_{int}$. Anything
satisfying these requirements is an instance of HYCO.

Everything else is an implementation choice adaptable to the problem at hand: the
ghost-point sampling used to approximate $\L_{int}$, the pretraining of the
synthetic model, the optimizer and learning-rate schedule, the stopping
criterion, and the specific values of $\alpha, \beta, H$. The choices made in our
experiments are recommended defaults, not part of the definition.

\subsection{Step-by-step recipe}
To apply HYCO to a new inverse or data-driven modeling problem, we suggest the
following steps.
\begin{enumerate}
    \item \textbf{Choose the physical model.} Any model producing a prediction
    from $\Lambda$ is admissible: a classical solver (finite differences, finite
    elements, finite volumes), a reduced-order model, or even a neural network or
    neural operator trained to minimize a PDE residual. A coarse, cheap solver
    often suffices; in the LWR example the physical mesh is forty times coarser
    than the one generating the data. The model must be differentiable with
    respect to its parameters $\Lambda$, since the $\Lambda$-update minimizes
    $\L_1$ by gradient descent; we obtain these gradients by automatic
    differentiation through the solver. A genuinely black-box solver can still be
    accommodated by supplying parameter gradients through finite differences or a
    gradient-free optimizer.
    \item \textbf{Choose the synthetic model.} Typically a neural network mapping
    the space-time coordinates to the solution. Unlike in PINNs, the network is
    never differentiated through the PDE operator, so its activation is
    unconstrained; we use ReLU networks throughout.
    \item \textbf{Allocate the data.} Decide which observations each model
    receives, which depends on the application. The synthetic model is always
    supplied with the data; the physical model may be given the data
    ($\beta > 0$) or driven by the interaction term alone ($\beta = 0$). In a
    privacy-aware setting the two models may even hold different, private
    datasets and still cooperate, since they exchange only predictions.
    \item \textbf{Define the interaction loss.} Choose a discrepancy between the
    two predictions (mean squared error in all our experiments) and a set of
    ghost points at which to evaluate it.
    \item \textbf{Set the ghost points.} Sample $H$ points in the domain and
    resample them at every epoch. For time-dependent problems, draw the temporal
    coordinates from the snapshot times stored by the solver, so that the
    physical prediction requires no temporal interpolation.
    \item \textbf{Alternate the updates.} Optionally pretrain the synthetic model
    on the data alone, then alternate the two minimizations in Gauss-Seidel
    fashion (Algorithm~\ref{alg:hyco}), one optimizer step per model per epoch.
    \item \textbf{Stop.} Terminate once the identified parameters stabilize,
    using the criterion \eqref{stopping}. An optional final refinement (e.g.
    L-BFGS-B with the physical parameters frozen) sharpens the synthetic fit.
    \item \textbf{Deploy.} Use the physical model for physically consistent
    extrapolation and the synthetic model for fast surrogate evaluation,
    depending on the application.
\end{enumerate}

\subsection{Choosing the hyperparameters}
The synthetic weight $\alpha$ and physical weight $\beta$ balance each model's
data term against the shared interaction loss, which carries unit weight in both.
In our experiments $\alpha \in [1, 10]$ is robust, with $\alpha = 10$ a good
default. Pushing $\alpha$ higher, as with $\alpha = 100$ in the Helmholtz sweep,
over-weights the synthetic data term: the network then interpolates the sensors
closely, but its poor extrapolation away from them is transmitted to the physical
model and degrades the recovered parameters. The weight $\beta$ should be set to
$0$ whenever the physical model receives no data, as in the Gray-Scott and LWR
experiments, where $\Lambda$ is driven by the interaction term alone. When the
physical model does see data, $\beta \in [0.1, 1]$ works well. Keeping $\beta$
moderate leaves the synthetic model responsible for interpolating the
observations, so that the physics is identified through agreement with a
synthetic prediction that already fits the data well, rather than by fitting the
sparse data directly.

The ghost-point count $H$ trades cost against the stability of the interaction
term. Final accuracy is insensitive to $H$ only up to a point: on Helmholtz even
$H = 50$ resampled points reach the same error floor, and $H = 200$ is a small,
safe default there. The number required grows with the complexity of the
solution, however; the richer Gray-Scott space-time field needs $H = 10^5$, so
more intricate solutions demand more ghost points. Beyond this, a larger $H$
mainly damps the oscillation of the interaction loss and stabilizes training,
at a cost that grows with $H$. Resampling the ghost points every epoch is
essential, as a fixed set causes premature lock-in of the parameters.

\subsection{When to use HYCO}
HYCO is at its best when a trustworthy physical model is available. Because the
two models regularize one another, a good physical ansatz pulls the synthetic
prediction toward physically consistent behavior; a poor or misspecified ansatz,
conversely, can drive the joint prediction away from the data, so the quality of
the physical model matters. 

Its advantage is most pronounced when the data is
localized: the interaction term then carries information from the data-rich
region, where the synthetic model interpolates well, into the unobserved region,
where the physical model extrapolates. This is precisely the regime of the LWR
experiment, whose sensors cover only part of the space-time domain. The framework
also handles non-smooth solutions that defeat smooth-ansatz physics-informed
methods, since the physical model resolves shocks the surrogate alone cannot.

Finally, when the inverse problem is nonconvex, classical solver-based inversion
such as Gauss-Newton or EKI converges in a few tens of iterations, at a small
fraction of HYCO's cost, but plateaus at a less accurate parameter set; the co-trained surrogate supplies a regularizing signal that a
purely local linearization lacks, and in our experiments this gap is not closed
by additional iterations. Conversely, when no fast solver is available the
physical update dominates the cost, and when the data is dense and the solution
smooth a single-model approach may already suffice.

\end{document}